\documentclass{article}

\usepackage{amssymb}
\usepackage{amsmath}
\usepackage[all]{xypic}
\usepackage[frenchb]{babel}
\usepackage{a4wide}
\usepackage{enumerate,fancyhdr}

\begin{document}

\title{Borne sur la torsion dans les vari\'et\'es ab\'eliennes de type CM}

\author{Nicolas Ratazzi \footnote{nicolas.ratazzi@math.u-psud.fr}
}

\renewcommand{\phi}{\varphi}
\renewcommand{\epsilon}{\varepsilon}
\renewcommand{\tilde}{\widetilde}
\newcommand{\Gal}{\textnormal{Gal}}
\newcommand{\tors}{\textnormal{tors}}
\newcommand{\Z}{\mathbb{Z}}
\newcommand{\Q}{\mathbb{Q}}
\newcommand{\C}{\mathbb{C}}
\newcommand{\R}{\mathbb{R}}
\newcommand{\N}{\mathbb{N}}
\renewcommand{\P}{\mathbb{P}}
\newcommand{\F}{\mathbb{F}}
\newcommand{\G}{\mathbb{G}}
\renewcommand{\hat}{\widehat}
\newcommand{\T}{\mathcal{T}}
\newcommand{\K}{\overline{K}}
\newcommand{\mP}{\mathcal{P}}
\newcommand{\mA}{\mathcal{A}}
\renewcommand{\O}{\mathcal{O}}
\renewcommand{\L}{\mathcal{L}}
\newcommand{\kk}{\mathbb{K}}
\newcommand{\im}{\textnormal{Im}}
\renewcommand{\H}{\textnormal{H}}
\newcommand{\MT}{\textnormal{MT}}
\newcommand{\Hom}{\textnormal{Hom}}
\newcommand{\Aut}{\textnormal{Aut}}
\newcommand{\hdg}{\textnormal{Hdg}}
\newcommand{\GL}{\textnormal{GL}}
\newcommand{\Res}{\textnormal{Res}}
\renewcommand{\ss}{\sigma}
\renewcommand{\text}{\textnormal}
\newcommand{\rrho}{\rho_{|\ell}}
\newcommand{\B}{\mathcal{B}}
\newcommand{\pr}{\textnormal{pr}}

\newtheorem{theo}{Th{\'e}or{\`e}me} [section]
\newtheorem{lemme}[theo]{Lemme}
\newtheorem{conj}{Conjecture}[section]
%\newcounter{nex}[theo]
\newtheorem{prop}[theo]{Proposition}
\newtheorem{cor}[theo]{Corollaire}

\newcommand{\defi}{\addtocounter{theo}{1}{\noindent \textbf{D{\'e}finition \thetheo\ }}}
\newcommand{\rem}{\addtocounter{theo}{1}{\noindent \textbf{Remarque \thetheo\ }}}

\newcommand{\demo}{\noindent \textit{D{\'e}monstration} : }

\maketitle

\hrulefill

\bigskip

\noindent \textbf{Abstract :} Let $A$ be an abelian variety of dimension $g\geq 1$ defined over a number field $K$. We study the size of the torsion group $A(F)_{\text{tors}}$ where $F/K$ is a finite extension and more precisely we study the best possible exponent $\gamma$ in the inequality Card$(A(F)_{\text{tors}})\ll [F:K]^{\gamma}$ when $F$ is any finite extension of $K$. In the CM case we give an exact formula for the exponent $\gamma$ in terms of the characters of the Mumford-Tate group---a torus in this case---and discuss briefly the general case.

\medskip

\noindent Finally we give an application of the main result in direction of a generalisation of the Manin-Mumford conjecture.

\bigskip

\hrulefill

\bigskip

\noindent \textbf{R\'esum\'e : } Soit $A$ une vari\'et\'e ab\'elienne de dimension $g\geq 1$ d\'efinie sur un corps de nombres $K$. On \'etudie la taille du groupe de torsion $A(F)_{\text{tors}}$ o\`u $F/K$ est une extension finie, et on \'etudie plus pr\'ecis\'ement le meilleur exposant possible $\gamma$ dans l'in\'egalit\'e  Card$(A(F)_{\text{tors}})\ll [F:K]^{\gamma}$ quand $F$ parcourt les extensions finies de $K$. Dans le cas CM, nous donnons une formule exacte pour l'exposant $\gamma$ en fonction des caract\`eres du groupe de Mumford-Tate---un tore dans ce cas---et nous donnons une br\`eve discussion dans le cas g\'en\'eral.

\medskip

\noindent Enfin nous donnons une application du r\'esultat principal en direction d'une g\'en\'eralisation de la conjecture de Manin-Mumford.

\bigskip

\hrulefill

\bigskip

\noindent \textit{classification :} 11G10, 11G15, 14K15, 11F80

\bigskip

\section{Introduction et r\'esultats}

Soient $K$ un corps de nombres et $A/K$ une vari\'et\'e ab\'elienne sur $K$ de dimension $g\geq 1$. Le classique th\'eor\`eme de Mordell-Weil assure que le groupe $A(K)$ des points $K$-rationnels de $A$ est de type fini. Un probl\`eme naturel qui se pose alors est de comprendre le sous-groupe de torsion $A(K)_{\tors}$. Dans le cas o\`u $A$ est une courbe elliptique (d\'efinie) sur $\Q$, Mazur \cite{mazur} a classifi\'e les groupes de torsion possibles. Ceci \'etant, ce probl\`eme semble tout \`a fait hors de port\'ee dans le cas g\'en\'eral et un sous-probl\`eme plus raisonnable consiste \`a essayer de comprendre le cardinal de $A(K)_{\tors}$ lorsque $A$ et $K$ varient. Il y a essentiellement deux approches possibles pour ce probl\`eme : soit l'on fixe le corps de nombres $K$ et l'on s'int\'eresse \`a la variation du cardinal lorsque $A$ d\'ecrit les vari\'et\'es ab\'eliennes sur $K$ de dimension $g$ fix\'ee. Dans cette direction citons le c\'el\`ebre r\'esultat de Merel \cite{merel} : si $E$ est une courbe elliptique sur un corps de nombres $K$, le cardinal de $E(K)_{\tors}$ est born\'e par une constante ne d\'ependant que du degr\'e de $K$ sur $\Q$. Parent \cite{parent} a rendu effectif le r\'esultat de Merel, obtenant une borne doublement exponentielle en le degr\'e $[K:\Q]$. En dimension sup\'erieure essentiellement rien n'est connu concernant ce probl\`eme connu sous le nom de conjecture de borne uniforme. La seconde approche possible concernant le cardinal de $A(K)_{\tors}$ consiste \`a fixer une vari\'et\'e ab\'elienne $A$ d\'efinie sur un corps de nombres $K_0$ et \`a faire varier $K$ parmi les extensions finies de $K_0$ ; l'objectif \'etant cette fois-ci d'obtenir une borne par une constante $C(A/K_0,[K:K_0])$ avec une d\'ependance explicite (la meilleure possible) en le degr\'e $[K:K_0]$ (ou, ce qui revient au m\^eme en le degr\'e $[K:\Q]$). Concernant ce second probl\`eme, Masser \cite{lettre} et \cite{mas} a montr\'e dans le cas g\'en\'eral que la d\'ependance est polynomiale en $[K:\Q]$. La question naturelle qui se pose est alors de savoir quel est le plus petit exposant $\gamma(A)$ possible dans cette borne polynomiale. Nous donnons dans cet article une r\'eponse \`a cette question dans le cas des vari\'et\'es ab\'eliennes CM. De plus ind\'ependamment de son int\'er\^et propre, l'obtention d'une borne meilleure que celle de Masser a des cons\'equences concr\`etes concernant des probl\`emes de g\'eom\'etrie diophantienne (\textit{cf.} th\'eor\`eme \ref{ray} ci-dessous).

\medskip

\noindent Soit $A/K$ une vari\'et\'e ab\'elienne de dimension $g\geq 1$ sur un corps de nombres $K$. On utilise la notation $\ll$ pour dire \`a une constante pr\`es ne d\'ependant que de $A/K$ et on pose 
\[\gamma(A)=\inf\left\lbrace x>0\, | \,  \forall F/K \text{ finie, }\ \left|A(F)_{\tors}\right|\ll [F:K]^x\right\rbrace.\]

\medskip

\noindent Dans le cas g\'en\'eral, la meilleure estimation de ce nombre est due \`a Masser \cite{lettre} et \cite{mas}.

\medskip

\begin{theo}\textnormal{\textbf{(Masser)}}\label{tmas} On a $\gamma(A)\leq g$. 
\end{theo}

\medskip

\noindent Dans le cas des courbes elliptiques de type CM, le r\'esultat de Masser est optimal. On peut se demander ce qu'il en est en dimension sup\'erieure. Pour cela nous avons besoin de rappeler la notion de groupe de Mumford-Tate d'une vari\'et\'e ab\'elienne $A/K$ sur un corps de nombres. 

\medskip

\noindent Fixons d\'esormais un plongement $K\subset \C$ et notons $\overline{K}$ une cl\^oture alg\'ebrique de $K$ dans $\C$. Soit $A/K$ une vari\'et\'e ab\'elienne. On note $V=H^1(A(\C),\Q)$ le premier groupe de cohomologie singuli\`ere de la vari\'et\'e analytique complexe $A(\C)$. C'est un $\Q$-espace vectoriel de dimension $2g$. Il est naturellement muni d'une structure de Hodge de type $\{(1,0),(0,1)\}$, c'est-\`a-dire d'une d\'ecomposition sur $\C$ de $V_{\C}:=V\otimes_{\Q}\C$ donn\'ee par $V_{\C}=V^{1,0}\oplus V^{0,1}$ telle que $V^{0,1}=\overline{V^{1,0}}$ o\`u $\overline{\hspace{.1cm} \cdot \hspace{.1cm}}$ d\'esigne la conjugaison complexe. On note $\mu : \G_{m,\C}\rightarrow \GL_{V_{\C}}$ le cocaract\`ere tel que pour tout $z\in\C^{\times}$, $\mu(z)$ agit par multiplication par $z$ sur $V^{1,0}$ et agit trivialement sur $V^{0,1}$. On d\'efinit le groupe de Mumford-Tate en suivant \cite{pink1}.
\medskip

\defi Le \textit{groupe de Mumford-Tate} $\MT(A)/\Q$ de $A$ est le plus petit $\Q$-sous-groupe alg\'ebrique $G$ de $\GL_V$ (vu comme $\Q$-sch\'ema en groupes) tel que, apr\`es extension des scalaires \`a $\C$, le cocaract\`ere $\mu$ se factorise \`a travers $G_{\C}:=G\times_{\Q}\C$.

\medskip

\noindent Nous faisons au paragraphe \ref{mtt} des rappels concernant le groupe de Mumford-Tate, dans le cas g\'en\'eral et plus sp\'ecifiquement dans le cas CM.

\medskip

\noindent Dans le cas des vari\'et\'es ab\'eliennes simples de type CM, Ribet a donn\'e une minoration de $\gamma(A)$. Pr\'ecis\'ement, en notant $\omega(n)$ le nombre de facteurs premiers de l'entier $n$ et $A[n]$ l'ensemble des points de $A(\overline{K})$ d'ordre divisant $n$, Ribet \cite{ribet} montre que 

\begin{theo}\textnormal{\textbf{(Ribet)}}\label{trib} Si $A/K$ est une vari\'et\'e ab\'elienne de type CM, alors, il existe deux constantes strictement positives $C_1$ et $C_2$  ne d\'ependant que de $A$ et $K$ telles que : pour tout entier $n\geq 1$, 
\[C_1^{\omega(n)}\leq\frac{[K(A[n]):K]}{n^{d}}\leq C_2^{\omega(n)},\]
\noindent o\`u $d=\dim \MT(A)$. De plus si $A$ est g\'eom\'etrique\-ment simple on a $d\geq 2+\log_2 g$.
\end{theo}

\medskip

\noindent Comme corollaire du th\'eor\`eme \ref{trib}, on obtient imm\'ediatement, dans le cas o\`u $A/K$ est une vari\'et\'e ab\'elienne de type CM, l'in\'egalit\'e 
\begin{equation}\label{tribcor}
\gamma(A)\geq \frac{2\dim A}{\dim \MT(A)}
\end{equation}
\noindent Nous montrons plus g\'en\'eralement que cette minoration reste valable pour toute vari\'et\'e ab\'elienne sur un corps de nombres. C'est l'objet du th\'eor\`eme suivant, prouv\'e au paragraphe \ref{1}.

\medskip

\begin{theo}\label{bla}Soit $A/K$ une vari\'et\'e ab\'elienne quelconque de dimension $g$. On a 
\[\gamma(A)\geq \frac{2g}{d}\]
\noindent o\`u $d$ est la dimension du groupe de Mumford-Tate de $A$.
\end{theo}

\medskip

\defi On dit qu'une vari\'et\'e ab\'elienne $A/K$ est \textit{sans facteur carr\'e} si elle est isog\`ene sur $\overline{K}$ \`a un produit $\prod_{i=1}^dA_i$ avec les $A_i$ simples et deux \`a deux non isog\`enes.

\medskip

\noindent Rappelons la notion de groupe des caract\`eres d'un tore alg\'ebrique sur un corps $k$ (\textit{i.e.} d'un groupe alg\'ebrique $G/k$ isomorphe sur $\overline{k}$ au groupe multiplicatif $\G_{m,\overline{k}}^{\dim G}$).

\medskip

\defi Soit $\T/k$ un tore alg\'ebrique sur un corps $k$. On note $\overline{k}$ la cl\^oture s\'eparable de $k$. On appelle \textit{groupe des caract\`eres de $\T$} et on note $X^*(\T)$ le groupe
\[X^*(\T)=\textnormal{Hom}_{\overline{k}}\left(\T_{\overline{k}}, \G_{m,\overline{k}}\right).\]
\noindent On d\'efinit de m\^eme le \textit{groupe des cocaract\`eres de $\T$} et on note $X_*(\T)$ le groupe
\[X_*(\T)=\textnormal{Hom}_{\overline{k}}\left(\G_{m,\overline{k}}, \T_{\overline{k}}\right).\]
\noindent Dans les deux cas pr\'ec\'edents et dans la suite, $\textnormal{Hom}_{\overline{k}}$ d\'esigne le groupes des homomorphismes entre groupes alg\'ebriques sur $\overline{k}$. On note $X^*(\T)\otimes\Q$ le $\Q$-espace vectoriel d\'eduit de $X^*(\T)$.

\medskip

\noindent Nous pouvons maintenant en venir au r\'esultat principal de cet article : l'obtention d'une valeur exacte pour $\gamma(A)$ dans le cas o\`u $A/K$ est une vari\'et\'e ab\'elienne de type CM, sans facteur carr\'e. De fait, le cas o\`u $A$ est g\'eom\'etriquement simple est d\'ej\`a suffisant dans les applications (\textit{cf.} par exemple le th\'eor\`eme \ref{ray} ci-dessous) et la proposition \ref{prop1} du paragraphe \ref{1} montre qu'il est facile d'obtenir un encadrement de $\gamma(A)$ en fonction des facteurs $\gamma(A_i)$ correspondant aux facteurs g\'eom\'etriquement simples $A_i$ de $A$. N\'eanmoins la preuve \'etant valable sans complication suppl\'ementaire dans le cas (l\'eg\`erement) plus g\'en\'eral des vari\'et\'es ab\'eliennes de type CM sans facteur carr\'e nous \'enon\c{c}ons notre r\'esultat dans ce cadre.

\medskip

\noindent Si $A/K$ est de type CM sans facteur carr\'e, on note $\T$ son groupe de Mumford-Tate~: c'est un tore alg\'ebrique sur $\Q$ (\textit{cf.} le sous-paragraphe \ref{rappels}). On note 
\[I=\left\lbrace\chi_1,\ldots,\chi_{2g}\right\rbrace\] 
\noindent l'ensemble des caract\`eres diagonalisant (sur $\overline{\Q}$) l'action de $\T$ sur l'espace vectoriel $V_{\overline{\Q}}$ de dimension $2g$ donn\'e par $V=H^1(A(\C),\Q)$ . L'hypoth\`ese \textit{sans facteur carr\'e} est faite pour assurer que les $\chi_i$ sont deux \`a deux distincts. 

\medskip

\defi Soit $A/K$ de type CM sans facteur carr\'e. Soit $W$ un sous-es\-pace vectoriel sur $\Q$ de $X^*(\T)\otimes \Q$. On pose
\[n(W)=\left|I\cap W\right|=\text{Card}\left(\left\lbrace i\in\{1,\ldots, 2g\} | \, \chi_i\in W  \right\rbrace\right).\]
\defi On d\'efinit un invariant $\alpha(A)$ comme suit :
\[\alpha(A)=\sup\left\lbrace \frac{n(W)}{\dim W}\, | \, W \text{ sous-$\Q$-espace vectoriel non nul de }X^*(\T)\otimes\Q \right\rbrace.\]

\rem On voit sur la d\'efinition que le nombre $\alpha(A)$ ne d\'epend que des caract\`eres $\chi_1,\ldots,\chi_{2g}$. En particulier, \'etant donn\'ee une vari\'et\'e ab\'elienne de type CM, \textit{i.e.} \'etant donn\'e un type CM le nombre $\alpha(A)$ est calculable explicitement et algorithmiquement. De fait il est clair que si l'on note $W_S=\text{Vect}(S)$ pour tout sous-ensemble $S$ de $I$, alors 
\[\alpha(A)=\max\frac{n(W_S)}{\dim W_S}\]
\noindent le max portant sur la collection finie des sous-ensembles $S$ de $I$. Ceci se calcule en d\'etermi\-nant les relations entre les diff\'erents caract\`eres $\chi\in I$.

\medskip

\noindent Avec ces notations, suivant une strat\'egie sugg\'er\'ee par Serre, notre r\'esultat principal est le suivant :

\medskip

\begin{theo}\label{thp}Soit $A/K$ une vari\'et\'e ab\'elienne sans facteur carr\'e, de type CM. On a 
\[\gamma(A)=\alpha(A).\]
\end{theo}

\medskip

\rem Notons quelque chose qui n'\'etait pas \'evident \textit{a priori} sur la d\'efinition de $\gamma(A)$~: dans le cas de type CM sans facteur carr\'e, l'exposant $\gamma(A)$ est un nombre rationnel (puisqu'il est clair sur la d\'efinition que $\alpha(A)$ l'est, avec un num\'erateur compris entre $2$ et $2g$ et un d\'enominateur compris entre $1$ et $\dim \T$).

\medskip

\noindent Par ailleurs ce r\'esultat est d'autant plus int\'eressant qu'il est possible de calculer une bonne majoration de $\alpha(A)$ en fonction de $g$ voire dans certains cas particuliers de calculer sa valeur exacte en fonction de $g$ et $d=\dim \T$.

\medskip

\begin{theo}\label{th2}Soit $A/K$ une vari\'et\'e ab\'elienne sans facteur carr\'e, de type CM, de dimension $g\geq 1$. On a
\[\alpha(A)\leq \frac{2g}{2+\log_2(g)},\]
\noindent o\`u l'on a not\'e $\log_2$ le logarithme en base $2$.
\end{theo}

\medskip

\noindent En cons\'equence de nos th\'eor\`emes \ref{thp} et \ref{th2}, nous obtenons :

\medskip

\begin{cor}\label{resco} Soit $A/K$ une vari\'et\'e ab\'elienne de type CM de dimension $g\geq 1$, sans facteur carr\'e. On a 
\[\gamma(A)\leq \frac{2g}{2+\log_2(g)},\]
\noindent o\`u l'on a not\'e $\log_2$ le logarithme en base $2$.
\end{cor}

\medskip

\rem On sait (\textit{cf.} par exemple \cite{dodson} theorem 1.0) que pour tout $g$ de la forme $2^n$, avec $n\geq 2$, il existe une vari\'et\'e ab\'elienne CM simple telle que la dimension de son groupe de Mumford-Tate est pr\'ecis\'ement $2+\log_2(g)$. En utilisant le th\'eor\`eme \ref{bla}, ceci prouve que la borne du corollaire pr\'ec\'edent sur $\gamma(A)$ est optimale en g\'en\'eral.

\medskip

\noindent On voit ainsi que d\`es lors que la dimension de la vari\'et\'e ab\'elienne est strictement sup\'erieure \`a $1$, ceci raffine le r\'esultat de Masser (dans le cas de type CM)~:

\medskip

\begin{cor} Soit $A/K$ une vari\'et\'e ab\'elienne sans facteur carr\'e, de type CM et de dimension $g\geq 2$. Alors, 

\[\gamma(A)< g.\]
\end{cor}

\medskip

\noindent Passons maintenant aux cas particuliers dans lesquels on peut calculer explicitement la constante $\alpha(A)$ et donc $\gamma(A)$. On introduit pour cela une d\'efinition :

\medskip

\defi On dit qu'une vari\'et\'e ab\'elienne $A$ de type CM, de dimension $g$ est de type \textit{non d\'eg\'en\'er\'e} si la dimension du groupe de Mumford-Tate $\T$ de $A$ est $g+1$.

\medskip

\rem \label{cds} Soit $A$ une vari\'et\'e ab\'elienne quelconque de dimension $g$, de groupe de Mumford-Tate de dimension $d$. On sait que $d\leq g+1$. Par ailleurs si $A$ est g\'eom\'etriquement simple de type CM, alors Ribet \cite{ribet} a montr\'e que 
\[2+\log_2(g)\leq d.\]
\noindent De plus un autre r\'esultat de Ribet (\textit{cf.} \cite{ribet2} Theorem 2) dit que si $A$ est une vari\'et\'e ab\'elienne g\'eom\'etriquement simple de type CM et de dimension un entier $g$ premier, alors $A$ est de type non d\'eg\'en\'er\'e.

\medskip

\begin{prop}\label{parti1}Soit $A/K$ une vari\'et\'e ab\'elienne g\'eom\'etriquement simple, de type CM, de dimension $g$. Si le type de $A$ est non d\'eg\'en\'er\'e ou si $g$ est inf\'erieur ou \'egal \`a $7$, alors
\[\alpha(A)=\frac{2g}{d}\]
\noindent o\`u $d$ est la dimension du groupe de Mumford-Tate de $A$.
\end{prop}

\medskip

\noindent Enfin il y a \'egalement un dernier cas o\`u l'on sait calculer la valeur de $\gamma(A)$ : si $A$ est une courbe elliptique sans multiplication complexe. Dans ce cas on conna\^it le groupe de Mumford-Tate de $A$, c'est le groupe alg\'ebrique $\text{GL}_2$ sur $\Q$. 

\medskip

\begin{prop}\label{parti2}Si $E/K$ est une courbe elliptique sans multiplication complexe, alors 
\[\gamma(E)=\frac{1}{2}.\]
\end{prop}

\medskip

\noindent Remarquons que dans ce dernier cas on a $\frac{2g}{d}=\frac{1}{2}$ et on trouve encore l'\'egalit\'e $\gamma(A)=\frac{2g}{d}$. Ainsi au vu des deux propositions \ref{parti1} et \ref{parti2} pr\'ec\'edentes et au vu de notre th\'eor\`eme \ref{bla}, on est tent\'e de poser la question suivante :

\medskip

\noindent \textbf{Question 1 : }si $A/K$ est une vari\'et\'e ab\'elienne g\'eom\'etriquement simple et $n\geq 1$, a-t-on 
\[\gamma(A^n)=\frac{2n\dim A}{\dim \MT(A)}\ ?\]

\noindent En plus des r\'esultats pr\'ec\'edents, cet \'enonc\'e se ram\`ene ais\'ement au cas o\`u $n=1$~: on v\'erifie dans le lemme \ref{puissance0} et le corollaire \ref{puissance1} que 
\[\forall n\geq 1\ \ \ \ \gamma(A^n)=n\gamma(A)\ \ \text{ et }\ \ \ \MT(A^n)\simeq\MT(A).\]
\noindent Ceci \'etant, un travail en cours semble sugg\'erer que cet \'enonc\'e sera faux en g\'en\'eral (pour une vari\'et\'e ab\'elienne quelconque) ; le bon \'enonc\'e dans le cas g\'en\'eral semblant plut\^ot \^etre le suivant : 

\medskip

\noindent \textbf{Question 2 : }si $A/K$ est une vari\'et\'e ab\'elienne isog\`ene au produit $\prod_{i=1}^rA_i^{n_i}$ o\`u les $A_i$ sont g\'eom\'etriquement simples et deux \`a deux non isog\`enes sur $\overline{\Q}$, a-t-on 
\[\gamma(A)=\max_{\emptyset\not=I\subset\{1,\ldots,r\}}\frac{2\dim \prod_{i\in I} A_i^{n_i}}{\dim \MT(\prod_{i\in I}A_i)}\ ?\]
\noindent Si $A$ est de la forme $A_1^n$ on retombe sur la premi\`ere question. Notons que d\'ej\`a dans le cas des vari\'et\'es ab\'eliennes de type CM, il serait tr\`es int\'eressant de savoir r\'epondre \`a ces questions : affirmativement ou en donnant un contre-exemple. Un autre test int\'eressant serait le cas d'un produit $E_1^{n_1}\times E_2^{n_2}$ avec $E_1,E_2$ deux courbes elliptiques quelconques et $n_1, n_2$ deux entiers.

\medskip

\noindent Pour conclure cette introduction nous donnons un exemple d'application des th\'eor\`emes \ref{thp} et \ref{th2} concernant une r\'ecente conjecture g\'en\'eralisant les conjectures de Manin-Mumford et de Mordell-Lang~:

\medskip

\noindent Soient $A/\overline{\Q}$ une vari\'et\'e ab\'elienne (d\'efinie) sur $\overline{\Q}$, $X/\overline{\Q}$ une courbe irr\'educ\-tible de $A$ et $r$ un entier positif ou nul. Suivant Bombieri, Masser et Zannier \cite{BMZ} dans le cas de $\mathbb{G}_m^n$ et R\'emond \cite{remond} dans le cas des vari\'et\'es ab\'eliennes, on s'int\'eresse au probl\`eme suivant : on consid\`ere l'ensemble
\[A^{[r]}:=\bigcup_{\text{codim\,}G\geq r}G(\overline{\Q})\]
\noindent o\`u l'union porte sur les sous-groupes alg\'ebriques non n\'ecessairement connexes de $A$ de codimension au moins $r$. \`A quelles conditions sur $r$ et sur $X$ peut-on garantir que l'ensemble $X(\overline{\Q})\cap A^{[r]}$ est fini ? C'est essentiellement \`a ce probl\`eme qu'est consacr\'e l'article de R\'emond. Notons que dans le cas le plus faible possible, si $r=\dim A$, on retombe sur un probl\`eme du type Manin-Mumford (th\'eor\`eme de Raynaud \cite{mray}). Si $A$ est une puissance d'une courbe elliptique, on peut voir que $X(\overline{\Q})\cap A^{[1]}$ est infini, donc on doit n\'ecessairement prendre $r\geq 2$. De mani\`ere ind\'ependante, Zilber (\cite{zilber} Conjecture 2) pour les vari\'et\'es semi-ab\'eliennes et Pink (\cite{pink} Conjecture 1.3) pour les vari\'et\'es de Shimura mixtes, ont formul\'e une conjecture qui, dans le cas des courbes incluses dans une vari\'et\'e ab\'elienne sur $\overline{\Q}$ se sp\'ecialise en la suivante : 

\medskip

\begin{conj}\label{condi}\textnormal{\textbf{(Zilber-Pink, cas particulier)}} Soient $A/\overline{\Q}$ une vari\'et\'e ab\'elienne et $X/\overline{\Q}$ une courbe dans $A$ irr\'eductible. Si $X$ n'est pas contenue dans un sous-groupe alg\'ebrique strict de $A$, alors l'ensemble $X(\overline{\Q})\cap A^{[2]}$ est fini.
\end{conj}

\medskip

\noindent R\'emond \cite{remond} a montr\'e que la conjecture pr\'ec\'edente est vraie si une tr\`es bonne minoration (conjecturale) des points d'ordre infini de $A$ est vraie : il s'agit d'une conjecture de David g\'en\'eralisant le probl\`eme de Lehmer. Par ailleurs dans le cas des vari\'et\'es ab\'eliennes de type CM, R\'emond \cite{remond} obtient un r\'esultat inconditionnel, mais sensiblement plus faible que la conjecture \ref{condi}~: soit $A$ une vari\'et\'e ab\'elienne de type CM, isog\`ene au produit $\prod_{i=1}^mA_i^{n_i}$ o\`u les $A_i$ sont des vari\'et\'es ab\'eliennes simples de dimension respective $g_i$, deux \`a deux non isog\`enes.

\medskip

\begin{theo}\textnormal{\textbf{(R\'emond) \cite{remond}}} Soient $A/\overline{\Q}$ une vari\'et\'e ab\'elienne de type CM et $X/\overline{\Q}$ une courbe dans $A$ qui n'est pas incluse dans un translat\'e de sous-vari\'et\'e ab\'elienne stricte de $A$, alors $X(\overline{\Q})\cap A^{[2+\sum_{i=1}^m g_i]}$ est fini.
\end{theo}

\medskip

\noindent Comme nous l'expliquons dans \cite{semilehmer}, en suivant la strat\'egie de R\'emond (et Bombieri-Masser-Zannier), en utilisant notre corollaire \ref{resco} ci-dessus concernant la borne sur les points de torsion pour les vari\'et\'es ab\'eliennes de type CM, ainsi qu'un r\'esultat de minoration de hauteur faisant l'objet de l'article s\'epar\'e \cite{semilehmer}, nous am\'eliorons ce th\'eor\`eme et obtenons un r\'esultat optimal dans le cas d'une puissance d'une vari\'et\'e ab\'elienne simple de type CM :

\medskip

\begin{theo}\label{ray} \textnormal{\textbf{(\cite{semilehmer})}} La conjecture \ref{condi} est vraie si $A$ est une puissance d'une vari\'et\'e ab\'elienne de type CM, simple de dimension quelconque.
\end{theo}

\medskip

\noindent Ceci g\'en\'eralise au cas d'une puissance d'une vari\'et\'e ab\'elienne de type CM simple de dimension quelconque le r\'esultat de Viada \cite{viada} et R\'emond-Viada (\cite{rv} th\'eor\`eme 1.7), valable pour une puissance d'une courbe elliptique \`a multiplication complexe. On peut m\^eme donner un r\'esultat un peu plus g\'en\'eral en fonction des exposants $\gamma(A_i)$ correspondant aux diff\'erents facteurs simples de $A$ : voir pour cela la remarque 1.5. de l'article \cite{semilehmer}.

\medskip

\noindent \textbf{Plan de l'article : }nous faisons au paragraphe \ref{mtt} les rappels n\'ecessaires concernant le groupe de Mumford-Tate $\MT(A)$, dans le cas g\'en\'eral et plus sp\'ecifiquement dans le cas CM. Nous donnons \'egalement le lien entre la repr\'e\-sen\-ta\-tion naturelle $\rho : \MT(A)\rightarrow \GL_V$ et les repr\'esentations $\ell$-adiques correspondant \`a l'action de $\Gal(\overline{K}/K)$ sur les points de torsion d'ordre une puissance de $\ell$ lorsque $\ell$ d\'ecrit l'ensemble des nombres premiers. Nous donnons au paragraphe \ref{1} une preuve de la minoration de $\gamma(A)$ annonc\'ee dans le th\'eor\`eme \ref{bla}. Ensuite nous donnons au paragraphe \ref{p2} une preuve du th\'eor\`eme \ref{th2}. Les paragraphes \ref{derp}, \ref{p3} et \ref{dernier} sont consacr\'es \`a la preuve du r\'esultat principal. Nous prouvons tout d'abord au paragraphe \ref{derp} l'in\'egalit\'e $\gamma(A)\geq \alpha(A)$ : la preuve est plus simple et contient d\'ej\`a une grande partie des id\'ees utilis\'ees dans la preuve de la seconde in\'egalit\'e. Nous discutons plus avant de la strat\'egie concernant la seconde in\'egalit\'e \`a la fin du paragraphe \ref{derp}. Enfin on s'int\'eresse aux cas particuliers correspondant aux propositions \ref{parti1} et \ref{parti2} dans le dernier paragraphe.

\medskip

\noindent \textbf{Remerciements : } Je tiens ici \`a exprimer ma sinc\`ere gratitude envers Jean-Pierre Serre. La preuve du th\'eor\`eme \ref{thp} est en effet bas\'ee sur une strat\'egie qu'il a eu la gentillesse de m'expliquer. Je tiens \'egalement \`a remercier le rapporteur pour son travail qui a permis d'am\'eliorer consid\'erablement la r\'edaction de cet article. Enfin je tiens \`a  remercier Marc Hindry pour l'aide pr\'ecieuse qu'il m'a apport\'ee : sans lui cet article n'aurait sans doute pas \'et\'e achev\'e.

\section{Rappels sur le groupe de Mumford-Tate et sur les tores alg\'ebriques\label{mtt}}
\subsection{Groupe et conjecture de Mumford-Tate}

\noindent  Concernant les groupes de Mumford-Tate $\MT(A)$ la d\'efinition a \'et\'e rappel\'ee dans l'introduction. 

\begin{lemme}\label{puissance0} Soit $A/K$ une vari\'et\'e ab\'elienne sur un corps de nombres $K$ plong\'e dans $\C$. Soit $n\geq 1$ un entier. On a 
\[\MT(A^n)\simeq\MT(A).\]
\end{lemme}
\demo On commence par remarquer que $H^1(A^n(\C),\Q)\simeq H^1(A(\C),\Q)^{\oplus n}$. Posons $V=H^1(A(\C),\Q)$. Vu la d\'efinition du groupe de Mumford-Tate, on voit que l'on obtient le r\'esultat en plongeant diagonalement $\GL_V$ dans $\GL_{V^{\oplus n}}$.\hfill$\Box$

\medskip

\noindent Rappelons maintenant la c\'el\`ebre conjecture de Mumford-Tate. Notons $A/K$ une vari\'et\'e ab\'elienne quelconque de dimension $g\geq 1$, d\'efinie sur un corps de nombres $K$ que l'on suppose plong\'e dans $\C$. Notons \'egalement $\MT(A)$ le groupe de Mumford-Tate correspondant. 

\medskip

\defi\label{rol} Notons $T_{\ell}(A)$ le module de Tate et $V_{\ell}:=T_{\ell}(A)\otimes_{\Z_{\ell}}\Q_{\ell}$. Soit $\ell$ un premier et $\rho_{\ell} : \Gal(\overline{K}/K)\rightarrow \Aut(T_{\ell}(A))\subset\GL(V_{\ell})$ la repr\'esen\-ta\-tion $\ell$-adique associ\'ee \`a l'action de Galois sur les points de $\ell^{\infty}$-torsion de $A$. On d\'efinit $\overline{G_{\ell}}$ comme \'etant l'adh\'erence de Zariski de l'image $G_{\ell}$ de $\rho_{\ell}$ dans le groupe alg\'ebrique $\GL_{V_{\ell}}\simeq \GL_{2g,\Q_{\ell}}$. C'est un groupe alg\'ebrique sur $\Q_l$ dont on notera $\overline{G_{\ell}}^0$ la composante neutre (composante connexe de l'identit\'e).

\medskip

\noindent Les th\'eor\`emes de comparaison entre cohomologie \'etale et cohomologie classique (\textit{cf.} \cite{sga4} XI) d'une part, et la comparaison entre le premier groupe de cohomologie \'etale et le module de Tate (\textit{cf.} \cite{milne} 15.1 (a)) d'autre part donnent pour tout premier $\ell$ l'isomorphisme canonique : 
\[V_{\ell}=T_{\ell}(A)\otimes_{\Z_{\ell}}\Q_{\ell}\simeq V\otimes_{\Q}\Q_{\ell}.\]
\noindent Nous fixons une fois pour toute dans la suite un tel isomorphisme. Ceci permet de comparer $\MT(A)\times_{\Q}\Q_{\ell}$ et $\overline{G_{\ell}}^0$~:

\begin{conj}\label{mt}\textnormal{\textbf{(Mumford-Tate)}} Pour tout $\ell$ premier, on a $\overline{G_{\ell}}^0=\MT(A)\times_{\Q}\Q_{\ell}$.
\end{conj}

\medskip

\defi \label{repmt} Notons $\rho : \MT(A)\hookrightarrow \GL_V$ la repr\'esentation naturelle du groupe de Mumford-Tate.

\medskip

\noindent Un certain nombre de cas particuliers, ainsi que de r\'esultats en direction de la conjecture pr\'ec\'edente sont connus. Nous renvoyons \`a la r\'ef\'erence \cite{pink1} pour une discussion d\'etaill\'ee de ces r\'esultats. Disons simplement ici que cette conjecture est un th\'eor\`eme pour les vari\'et\'es ab\'eliennes de type CM (travaux de Shimura-Taniyama \cite{shita}) ainsi que pour les courbes elliptiques sans multiplication complexe (Serre \cite{serrenoncm}). De mani\`ere g\'en\'erale on sait par les travaux de Borovo{\u\i} \cite{boro}, Deligne \cite{del} Exp I, 2.9, 2.11, et Pjatecki{\u\i}-{\v{S}}apiro \cite{piat} qu'une inclusion est toujours vraie~:

\medskip

\begin{theo}\label{bdp}Pour tout $\ell$ premier, $\overline{G_{\ell}}^0\subset \MT(A)\times_{\Q}\Q_{\ell}$.
\end{theo}

\medskip 

\noindent Par ailleurs, on a le r\'esultat suivant, d\^u \`a Shimura-Taniyama \cite{shita} dans le cas de type CM et d\^u \`a Serre \cite{college} 2.2.3. (\textit{cf.} \'egalement \cite{serken}), dans le cas g\'en\'eral.

\medskip

\begin{theo} L'application 
\[\epsilon~: \Gal(\overline{K}/K)\longrightarrow \overline{G_{\ell}}(\Q_{\ell})\longrightarrow\overline{G_{\ell}}(\Q_{\ell})/\overline{G_{\ell}}^0(\Q_{\ell})\]
\noindent est continue surjective, de noyau ind\'ependant de $\ell$ pour tout premier $\ell$.
\end{theo}

\medskip

\noindent Ainsi, le noyau de $\epsilon$ est un sous-groupe d'indice fini de $\Gal(\overline{K}/K)$ donc il existe une extension finie $K'/K$ telle que $\Gal(\overline{K'}/K')\subset \ker\epsilon$. Dit autrement, les deux th\'eor\`emes pr\'ec\'edents donnent le lien entre la repr\'esentation $\rho$ et les repr\'esentations $\ell$-adiques $\rho_{\ell}$ pour tout premier $\ell$ : quitte \`a remplacer au d\'epart $K$ par une extension $K'$ finie ne d\'ependant que de $A$ (ce que nous ferons dans la suite), on a pour tout premier $\ell$ la factorisation
\[\rho_{\ell} : \Gal(\overline{K}/K)\rightarrow \MT(A)(\Q_{\ell})\overset{\rho}{\rightarrow}\GL_V(\Q_{\ell})\simeq \GL (V_{\ell}).\]

\noindent Enfin dans deux cas particuliers d\'ej\`a mentionn\'es pr\'ec\'edemment, on peut encore pr\'eciser les choses~:

\medskip

\begin{theo}\label{trr}\textnormal{\textbf{(Serre)}} Soit $A/K$ une vari\'et\'e ab\'elienne de type CM ou une courbe elliptique sans multiplication complexe. Pour tout $\ell$ premier on a
\[G_{\ell}=\im(\rho_{\ell})\subset \MT(A)(\Z_{\ell}),\]
\noindent cette inclusion \'etant de conoyau fini, born\'e ind\'ependamment de $\ell$. 
\end{theo}

\demo Dans le cas des vari\'et\'es ab\'eliennes de type CM cela d\'ecoule du th\'eor\`eme 1 p.II-26 de \cite{serabl}. Dans le cas des courbes elliptiques sans multiplication complexe il s'agit du th\'eor\`eme 3 p. 299 de \cite{serrenoncm}.\hfill$\Box$

\subsection{Groupe de Mumford-Tate des vari\'et\'es ab\'eliennes de type CM\label{rappels}}

\noindent Int\'eressons-nous maintenant plus particuli\`erement au cas des vari\'et\'es ab\'eliennes de type CM. Nous suivons pour cela la r\'ef\'erence \cite{murty} p.108-110. Commen\c{c}ons par le cas d'une vari\'et\'e g\'eom\'etri\-que\-ment simple de type CM de dimension $g$. Soit $E$ le corps CM de la vari\'et\'e ab\'elienne $A$. Rappelons que l'on a not\'e dans l'introduction $V=H^1(A(\C),\Q)$ le premier groupe de cohomologie singuli\`ere de la vari\'et\'e analytique complexe $A(\C)$. C'est un $\Q$-espace vectoriel de dimension $2g$. Il est naturellement muni d'une structure de Hodge de type $\{(1,0),(0,1)\}$, c'est-\`a-dire d'une d\'ecomposition sur $\C$ de $V_{\C}:=V\otimes_{\Q}\C$ donn\'ee par $V_{\C}=V^{1,0}\oplus V^{0,1}$ telle que $V^{0,1}=\overline{V^{1,0}}$ o\`u $\overline{\hspace{.1cm} \cdot \hspace{.1cm}}$ d\'esigne la conjugaison complexe. Mieux, suivant \cite{murty} \'equation (3) p. 108 on a la d\'ecomposition 
\[V^{1,0}=\bigoplus_{\sigma\in\Phi}V^{1,0}_{\sigma}\]
\noindent o\`u $\Phi$ est un type CM, \textit{i.e.} un ensemble de plongements $\sigma : E\hookrightarrow \C$ tel que 
\[\Phi\cup\overline{\Phi}=\Hom(E,\C)\ \ \text{  et  }\ \ \Phi\cap\overline{\Phi}=\varnothing.\]
\noindent Notons $T_E=\Res_{E/\Q}\G_{m, E}$ le tore correspondant au corps CM $E$. Si $\sigma$ est l'un des plongements de $E$ dans $\overline{\Q}$, alors $\sigma$ s'\'etend en un morphisme $E\otimes_{\Q}\overline{\Q}\rightarrow \overline{\Q}$ et d\'efinit donc un caract\`ere $[\sigma]\in X^*(T_E)$. Une base du groupe des caract\`eres $X^*(T_E)$ est donn\'ee par la famille $\{[\sigma]\ |\ \sigma : E\hookrightarrow \C\}$. Par d\'efinition $\MT(A)$ est le plus petit sous-groupe alg\'ebrique sur $\Q$ de $\GL_V$ contenant l'image de
\[\mu=\sum_{i=1}^ge_i\]
\noindent les $e_i\in X_*(T_E)$ \'etant les cocaract\`eres duaux des $[\sigma_i]\in X^*(T_E)$ avec $\{\sigma_1,\ldots,\sigma_g\}=\Phi$. En effet, $\mu$ n'est rien d'autre que l'homomorphisme d\'efini par la m\^eme notation dans l'introduction et donn\'e sur les complexes par : $\mu_{\C}~:\G_{m,\C}\rightarrow \GL_{V_{\C}}$ tel que pour tout $z\in\C^{\times}$, $\mu(z)$ agit par multiplication par $z$ sur $V^{1,0}$ et agit trivialement sur $V^{0,1}$. Le groupe $\MT(A)$ est inclus dans $T_E$, c'est donc un sous-tore alg\'ebrique : on le note d\'esormais $\T$ pour simplifier les notations. De plus par d\'efinition $\mu$ est un cocaract\`ere de $\T$. On note $\mu'$ son cocaract\`ere conjugu\'e (donn\'e par $\mu'=\sum_{i=g+1}^{2g}e_i$). Le groupe des caract\`eres $X^*(\T)$ de $\T$ est un quotient de $X^*(T_E)$. On note $\chi_1,\ldots,\chi_{2g}$ les images dans $ X^*(\T)$ des caract\`eres $[\sigma_1],\ldots,[\sigma_{2g}]$ de $X^*(T_E)$.

\medskip

\noindent Ce qui pr\'ec\`ede s'applique de m\^eme dans le cas d'une vari\'et\'e ab\'elienne sans facteur carr\'e de type CM, \`a condition de remplacer le corps CM $E$ par le produit de corps CM $\prod_i E_i$ et en rempla\c{c}ant le tore $T_E$ par le tore $\prod_iT_{E_i}$. On suppose d\'esormais que $A$ est de type CM sans facteur carr\'e.

\medskip

\defi \'Etant donn\'e un tore $T$ sur un corps $k$ de cl\^oture s\'eparable $\overline{k}$, notons $\langle\ , \ \rangle$ l'accouplement donn\'e par
\[X^*(T)\times X_*(T)\rightarrow \mathbb{Z}, \ \ \ (x,y)\mapsto \langle x,y\rangle:=\deg (x\circ y).\]

\begin{prop} \label{lh1}En \'etendant les scalaires \`a $\overline{\Q}$, les images de $\sigma \mu$, $\sigma$  d\'ecrivant $\Gal(\overline{\Q}/\Q)$, engendrent $\T_{\overline{\Q}}$. Par ailleurs, pour tout entier $i$ compris entre $1$ et $2g$ on a 
\[ \langle \chi_i, \mu+\mu'\rangle=1\ \text{ et }\ \langle\chi_i,\mu\rangle\in\{0,1\}.\]
\end{prop}
\demo Par construction du groupe de Mumford-Tate (\textit{cf.} par exemple \cite{pink1} Fact 5.9) on sait que les images de $\sigma \mu$, $\sigma$  d\'ecrivant $\Gal(\overline{\Q}/\Q)$, engendrent $\T_{\overline{\Q}}$, ceci \'etant pr\'ecis\'ement la premi\`ere assertion de la proposition. La seconde d\'ecoule directement des d\'efinitions et de la discussion pr\'ec\'edant l'\'enonc\'e de la proposition.\hfill$\Box$

\subsection{Tores alg\'ebriques}

\noindent Nous rappelons ici quelques r\'esultats et notations sur les tores alg\'ebriques.

\medskip

\begin{prop}\label{o1} Soient $k$ un corps de cl\^oture s\'eparable $\overline{k}$ et $T_1$, $T_2$ deux tores alg\'ebriques (d\'efinis) sur $k$. Alors, $T_1$ et $T_2$ sont isog\`enes sur $k$ si et seulement si les $\Q$-espaces vectoriels de caract\`eres $X^*(T_1)\otimes \Q$ et $X^*(T_2)\otimes \Q$ sont isomorphes en tant que $\Gal(\overline{k}/k)$-modules.
\end{prop}
\demo C'est la proposition 1.3.2. de \cite{ono}.\hfill$\Box$

\medskip

\noindent \textbf{Notations }Soit $T$ un tore et $X^*(T)$ son groupe de caract\`eres. On se donne un sous-tore $T_1$ de $T$ et un sous-$\Z$-module libre $X_1$ de $X^*(T)$. On note alors
\[X_1^{\perp}=\bigcap_{\chi\in X_1}\ker \chi,\ \ \text{ et }\ \ T_1^{\perp}=\left\{\chi\in X^*(T)\, | \, T_1\subset \ker\chi\right\}.\]

\begin{prop}\label{o2}Avec les notations pr\'ec\'edentes, on a
\[X^*(T/T_1)\simeq T_1^{\perp}, \ \ \text{ et } X^*(T_1)\simeq X^*(T)/T_1^{\perp}.\]
\noindent De plus, on a 
\[(X_1^{\perp})^{\perp}= X_1,\ \ \text{ et }\ \ (T_1^{\perp})^{\perp}= T_1.\]
\end{prop}
\demo C'est la proposition 1.1.1. de \cite{ono}.\hfill$\Box$

\medskip

\noindent Le r\'esultat suivant est d\^u \`a Ribet. Nous introduisons pour cela une nouvelle notation :

\medskip

\noindent \textbf{Notation } Soient $\ell$ un nombre premier, $\kk$ une extension finie de $\Q$, contenue dans $\Q_{\ell}$, $X/\kk$ un tore alg\'ebrique et $\mathcal{X}$ un mod\`ele sur $\mathcal{O}_{\kk}$. Si $n$ est un entier strictement positif, on note 
\[X\left(\Z/\ell^n\Z\right):=\mathcal{X}(\Z_{\ell})/\mathcal{X}(1+\ell^n\Z_{\ell}).\]

\medskip

\begin{theo}\textnormal{\textbf{(Ribet)}}\label{rib} Soient $\ell$ un nombre premier, $\kk$ une extension finie de $\Q$, contenue dans $\Q_{\ell}$, $X/\kk$ un tore alg\'ebrique de dimension $\nu$ et $n$ un entier strictement positif. Il existe deux constantes $C$ et $C'$ strictement positives, ne d\'ependant que de $X/\kk$ telles que 
\[C' \ell^{n\nu}\geq X\left(\Z/\ell^n\Z\right)\geq C \ell^{n\nu}.\]
\end{theo}
\demo Il s'agit du th\'eor\`eme (2.5) de Ribet \cite{ribet}. Pr\'ecis\'ement, Ribet montre cet \'enonc\'e avec $\kk=\Q$ mais sa preuve vaut dans le cas \'enonc\'e ci-dessus.\hfill$\Box$

\section{Un encadrement et une minoration de $\gamma(A)$\label{1}}

\subsection{Un encadrement de $\gamma(A)$}

\noindent On donne ici un encadrement \'el\'ementaire de l'exposant $\gamma$ pour un produit de vari\'et\'es ab\'eliennes.

\medskip

\begin{prop}\label{prop1}Soient $r\geq 2$ un entier et $A_1,\ldots, A_r$ des vari\'et\'es ab\'eliennes sur un corps de nombres $K$. Soient $n_1,\ldots, n_r\geq 1$ des entiers. On a 
\[\max_{\emptyset\not=S\subset\{1,\ldots,r\}}\left(\min_{i\in S}n_i\right)\gamma\left(\prod_{j\in S}A_j\right)\leq\gamma\left(\prod_{i=1}^rA_i^{n_i}\right)\leq\sum_{i=1}^rn_i\gamma(A_i).\]
\end{prop}
\demo Cela d\'ecoule imm\'ediatement de ce que si $F/K$ est une extension finie et $P=(P_1,\ldots,P_m)$ est un point de $A_1\times\ldots\times A_m(F)$, alors $P$ est de torsion si et seulement si tous les $P_i$ le sont. Explicitons par exemple la preuve de l'in\'egalit\'e de gauche~: soit $S$ un sous-ensemble non vide de $\{1,\ldots,r\}$. On a
\begin{align*}
\left|\left(\prod_{j\in S}A_j\right)(F)_{\tors}\right|^{\min_{i\in S}n_i}	& = \left|\left(\prod_{j\in S}A_j^{\min_{i\in S}n_i}\right)(F)_{\tors}\right|\\
																			& \leq \left|\left(\prod_{i=1}^rA_i^{n_i}\right)(F)_{\tors}\right|.
\end{align*}
\noindent On en d\'eduit 
\[\gamma\left(\prod_{j\in S}A_j\right)\leq\frac{1}{\min_{i\in S}n_i}\gamma\left(\prod_{i=1}^rA_i^{n_i}\right).\]
\noindent Ceci prouve la premi\`ere in\'egalit\'e. La seconde se d\'emontre de la m\^eme fa\c{c}on.\hfill$\Box$

\medskip

\begin{cor}\label{puissance1}Soient $A/K$ une vari\'et\'e ab\'elienne sur un corps de nombres $K$ et $n\geq 1$ un entier. On a 
\[\gamma(A^n)=n\gamma(A).\]
\noindent De plus si $A$ est isog\`ene (sur $\overline{K}$) au produit $\prod_{i=1}^rA_i^{n_i}$ les $A_i$ \'etant deux \`a deux non isog\`enes, on a
\[\max\left\{\max_{1\leq i\leq r}n_i\gamma(A_i) , (\min_{1\leq i\leq r}n_i)\gamma\left(\prod_{i=1}^rA_i\right)\right\}\leq\gamma(A)\leq\sum_{i=1}^rn_i\gamma(A_i).\]
\end{cor}

\medskip

\noindent Ceci montre que pour avoir un encadrement de $\gamma(A)$ dans le cas g\'en\'eral, il suffit de savoir calculer $\gamma(A)$ pour les vari\'et\'es ab\'eliennes $A/K$ sans facteur carr\'e. Notons toutefois que l'encadrement donn\'e dans le corollaire peut tout de m\^eme \^etre tr\`es large. Ceci se voit par exemple en consid\'erant les deux vari\'et\'es ab\'eliennes $A_1=E_1\times E_0^{n_0}$ et $A_2=\prod_{i=1}^{n_1}E_i\times E_0^{n_0}$ o\`u les $E_i$ sont des courbes elliptiques CM deux \`a deux non isog\`enes d\'efinies sur un corps de nombres $K$ et o\`u $n_0$ et $n_1$ sont deux entiers strictement positifs tels que $n_0\geq 2(n_1+1)$. En admettant les r\'esultats de l'introduction (pr\'ecis\'ement le corollaire \ref{resco}), l'encadrement pr\'ec\'edent donne dans ces deux cas
\[n_0\leq\gamma(A_1) \leq n_0+1\ \text{ et}\ n_0\leq \gamma(A_2)\leq n_0+n_1.\]

\subsection{Minoration de $\gamma(A)$}

\noindent Comme annonc\'e on \'etend ici l'in\'egalit\'e (\ref{tribcor}) de l'introduction au cas d'une vari\'et\'e ab\'elienne quelconque.

\medskip

\begin{prop} Soit $A/K$ une vari\'et\'e ab\'elienne quelconque avec $K$ un corps de nombres plong\'e dans $\C$. On a
\[\gamma(A)\geq\frac{2\dim A}{\dim \MT(A)}\]
\noindent o\`u $\MT(A)$ est le groupe de Mumford-Tate de $A$.
\end{prop}
\demo On se donne un premier $\ell$, un entier strictement positif $n$ tels que $(\ell,n)\not=(2,1)$ et le groupe $H=A[\ell^n]$. On veut montrer que
\[\ell^{2gn}=|H|\geq c(A/K)[K(H):K]^{\frac{2g}{d}}\]
\noindent o\`u $g=\dim A$, $d=\dim \MT(A)$ et $c(A/K)$ est une constante ne d\'ependant que de $A/K$. Pour cela on introduit les voisinages de l'identit\'e suivants de $\text{GL}_{2g}(\Z_{\ell})$ : 
\[\forall n\in \N,\ \ \ V_n:=\left\{x\in \text{GL}_{2g}(\Z_{\ell})\ / \ x=\text{Id}\text{ mod }\ell^n\right\}, \ \ \ \text{ et }\ \ \ G_n=G_{\ell}\cap V_n.\]
\noindent On v\'erifie imm\'ediatement que pour tout entier $n$ positif, $G_n/G_{n+1}$ et $V_n/V_{n+1}$ sont des $\F_{\ell}$-espaces vectoriels, o\`u l'on a not\'e $\F_{\ell}$ le corps \`a $\ell$ \'el\'ements. Par ailleurs, on a les inclusions de groupes (c'est pour que la premi\`ere inclusion soit vraie que l'on a fait l'hypoth\`ese sur le couple $(\ell,n)$)~:
\[ G_n/G_{n+1}\hookrightarrow G_{n+1}/G_{n+2}, \ \ \ \text{ et }\ \ \ G_n/G_{n+1}\hookrightarrow V_n/V_{n+1},\] 
\noindent la premi\`ere inclusion \'etant d\'efinie par l'\'el\'evation \`a la puissance $\ell$ et la seconde \'etant induite par l'application $x\mapsto x$. De plus, $V_n/V_{n+1}$ s'injecte dans l'espace tangent de $\text{GL}_{2g}(\F_{\ell})$ en l'identit\'e. Notons 
\[d_n=\dim_{\F_{\ell}}G_n/G_{n+1}.\]
\noindent La suite $(d_n)_{n\in \N}$ est croissante et stationnaire \`a partir d'un certain rang $n_0$. On note $d_{\infty}$ sa limite. Le groupe $G_{n_0}$ est limite profinie des groupes $G_{n_0}/G_{n_0+k}.$ Dans ces conditions on peut montrer que le groupe $G_{n_0}$ est hom\'eomorphe \`a $\Z_{\ell}^{d_{\infty}}$. En particulier on en d\'eduit que 
\[d_{\infty}\leq \dim \overline{G_{\ell}}\leq\dim\MT(A)\]
\noindent la derni\`ere in\'egalit\'e r\'esultant du th\'eor\`eme \ref{bdp}. On a ainsi :
\[ [K(H):K]=\left|G_{\ell}/G_n\right|=\prod_{k=0}^{n-1}\left|G_i/G_{i+1}\right|\leq \ell^{nd_{\infty}}\leq \ell^{nd}.\]
\noindent On conclut en \'elevant ceci \`a la puissance $\frac{2g}{d}$.\hfill $\Box$

\section{Majoration de l'exposant $\alpha(A)$ : le th\'eor\`eme \ref{th2}\label{p2}}

\noindent Soient $K$ un corps de nombres et $A/K$ une vari\'et\'e ab\'elienne  de type CM, sans facteur carr\'e et de dimension $g$. Nous utiliserons dans la suite les notations du paragraphe \ref{rappels}. Notamment on note $\chi_1,\ldots,\chi_{2g}$ les caract\`eres de $\T$ diagonalisant l'action de $\T$ sur $V_{\overline{\Q}}$. L'hypoth\`ese faite sur $A$, \`a savoir qu'elle est sans facteur carr\'e, entra\^ine en particulier que les caract\`eres $\chi_i$ sont deux \`a deux distincts. Suivant une id\'ee de Ribet et Lenstra (\textit{cf.} \cite{ribet} p. 87), on peut en fait montrer mieux. Les caract\`eres $\chi_i$ sont dans $X^*(\T)$ qui est un $\Z$-module libre de rang fini. \'Etant donn\'e un nombre premier $\ell$, nous dirons que \textit{deux caract\`eres co\"{i}ncident modulo $\ell$} s'ils co\"{i}ncident dans le quotient $X^*(\T)/\ell X^*(\T)$.

\begin{lemme}\label{l1} Les caract\`eres $\chi_i$ sont deux \`a deux distincts modulo $2$.
\end{lemme}
\demo Soient $i$ et $j$ deux entiers tels que $\chi_i=\chi_j\mod 2$. Montrons que $\chi_i=\chi_j$. Par la proposition \ref{lh1} on sait que en \'etendant les scalaires \`a $\overline{\Q}$, les images de $\sigma \mu$, $\sigma$  d\'ecrivant $\Gal(\overline{\Q}/\Q)$, engendrent $\T_{\overline{\Q}}$.  Ainsi pour montrer que deux caract\`eres co\"{i}ncident, il suffit de montrer que 
\[\forall\sigma\in \Gal(\overline{\Q}/\Q),\ \ \ \langle\chi_i-\chi_j,\sigma \mu\rangle=0.\]
\noindent Or $\langle\chi_i-\chi_j,\sigma^{-1} \mu\rangle=\langle\chi_i^{\sigma}-\chi_j^{\sigma}, \mu\rangle$. De plus pour tout $\sigma$ et pour tout $\chi\in I$, $\chi^{\sigma}\in I$. D'autre part si deux caract\`eres $s,t$ co\"{i}ncident modulo $2$, alors $s^{\sigma}=t^{\sigma}\mod 2$. Il suffit donc de montrer que si $\chi_i=\chi_j\mod 2$ alors $\langle\chi_i-\chi_j, \mu\rangle=0$. Or par la proposition \ref{lh1} on sait que $\langle\chi_i,\mu\rangle\in\{0,1\}$, et de m\^eme pour $j$. On en d\'eduit donc que $\langle \chi_i-\chi_j,\mu\rangle\in\{-1,0,1\}$. On conclut en remarquant qu'un nombre dans cet ensemble vaut z\'ero modulo $2$ si et seulement s'il est nul. \hfill$\Box$

\medskip

\rem La m\^eme preuve montre plus g\'en\'eralement que les $\chi_i$ sont deux \`a deux distincts modulo $\ell$ pour tout nombre premier $\ell$.

\medskip

\begin{cor}\label{c1new}Soit $W$ un sous-$\Q$-espace vectoriel non nul de $X^*(\T)\otimes\Q$. En utilisant la notation $n(W)=\left|I\cap W\right|$ on a 
\[n(W)\leq 2^{\dim W-1}.\]
\end{cor}
\demo Quitte \`a renum\'eroter, soient $\chi_1,\ldots, \chi_{n(W)}$ les $n(W)$ \'el\'ements de $I$ dans $W$. Ils engendrent dans $X^*(\T)$ un $\Z$-module libre $X$ de rang inf\'erieur \`a $\dim W$. En r\'eduisant modulo $2$, on voit que $X/2X$ est un $\Z/2\Z$-espace vectoriel de dimension major\'ee par $\dim W$ et qui contient les $\bar \chi_1,\ldots,\bar\chi_{n(W)}$, o\`u $\bar\chi$ d\'esigne le caract\`ere $\chi$ modulo $2$. Par le lemme \ref{l1} pr\'ec\'edent les $\bar\chi_i$ sont deux \`a deux distincts. Par ailleurs, la proposition \ref{lh1} indique que si $\chi\in I$ alors $\langle\chi,\mu+\mu'\rangle=1$. Donc pour tout entier $i$ compris entre $1$ et $n(W)$, les $\chi_i$ sont contenus dans l'hyperplan  affine $H_1=\left\{\chi\, | \, \langle \chi, \mu+\mu'\rangle=1\right\}$ de $W$. En r\'eduisant modulo $2$, on conclut par cardinalit\'e.\hfill$\Box$

\medskip

\begin{cor}\label{c2}On a
\[\alpha(A)\leq \frac{2g}{2+\log_2g}\]
\noindent o\`u $\log_2$ d\'esigne le logarithme en base $2$.
\end{cor}
\demo Soit $W$ un sous-$\Q$-espace vectoriel non nul de $X^*(\T)\otimes\Q$. Il s'agit de montrer que $\frac{n(W)}{\dim W}\leq \frac{2g}{2+\log_2g}$. On d\'ecoupe la preuve en deux morceaux : 
\begin{enumerate}
\item Si $\dim W\leq 2+\log_2 g$, alors la croissance de la fonction $x\mapsto \frac{2^{x-1}}{x}$ permet de conclure.
\item Sinon on majore na\"ivement $n(W)$ par $2g$ ce qui permet encore de conclure.\hfill$\Box$
\end{enumerate}

\section{Le th\'eor\`eme \ref{thp}  : l'in\'egalit\'e $\alpha(A)\leq \gamma(A)$\label{derp}}

\noindent Avant de prouver l'in\'egalit\'e qui nous int\'eresse le plus, \`a savoir une majoration de l'exposant $\gamma(A)$, nous commen\c{c}ons par montrer dans ce paragraphe que la borne donn\'ee par $\alpha(A)$ est n\'ecessairement optimale : on a $\alpha(A)\leq \gamma(A)$. La preuve de ceci correspond \`a un cas particulier plus simple de la preuve de l'in\'egalit\'e r\'eciproque. Certaines id\'ees intervenant sous forme moins technique, il nous semble instructif de d\'etailler ici les choses, en esp\'erant rendre ainsi plus accessible la preuve de la seconde in\'egalit\'e.

\medskip

\noindent Rappelons bri\`evement les notations que nous utiliserons ici et dans la suite : $A/K$ est une vari\'et\'e ab\'elienne de dimension $g$, de type CM, sans facteur carr\'e sur un corps de nombres $K$. Son groupe de Mumford-Tate est not\'e $\T$ : c'est un tore alg\'ebrique sur $\Q$. Pour tout premier $\ell$ nous notons $V_{\ell}=T_{\ell}(A)\otimes \Q_{\ell}$ le module de Tate : c'est un $\Q_{\ell}$-espace vectoriel de dimension $2g$. Enfin nous notons, comme dans la d\'efinition \ref{rol}, $G_{\ell}$ l'image dans $\Aut(T_{\ell}(A))\simeq \GL_{2g}(\Z_{\ell})$  de la repr\'esentation $\ell$-adique $\rho_{\ell}$.

\medskip

\noindent Soit $W$ un sous-espace vectoriel non nul de $X^*(\T)\otimes \Q$ r\'ealisant le sup $\alpha(A)$. Soit $\ell$ un nombre premier totalement d\'ecompos\'e dans le produit de corps CM correspondant \`a $A$, de sorte que l'action de $\T$ sur $V_{\ell}$ est d\'ecompos\'ee : les caract\`eres $\chi_1,\ldots,\chi_{2g}$ diagonalisant l'action sont rationnels sur $\Q_{\ell}$. Quitte \`a renum\'eroter et \`a extraire une base, on a : 
\[W=\text{Vect}_{\Q}(\chi_1,\ldots,\chi_{\dim W}),\]
\noindent les $\chi_1,\ldots,\chi_{\dim W}$ formant une base de $W$. L'espace vectoriel $W$ \'etant fix\'e dans la suite nous noterons \'egalement
\[H=\left\{Q\in A[\ell]\, | \, Q=\sum_{i\in I_W}m_iP_i, \ m_i\in\F_{\ell}\right\},\]
\noindent o\`u $I_W=\{i\in [\![1,2g]\!]\, | \, \chi_i\in W\}$ est de cardinal $n(W)$ et $\{P_1,\ldots, P_{2g}\}$ est la base du $\F_{\ell}$-espace vectoriel $A[\ell]$ de dimension $2g$, dans laquelle la repr\'esentation se diagonalise.

\medskip

\noindent Posons $L=K(H)$. On a par d\'efinition de $H$~:
\[ |A(L)_{\tors}|\geq |H|=\ell^{n(W)}.\] 
\noindent Nous voulons montrer que $\ell^{n(W)}\gg [L:K]^{\alpha(A)}=[L:K]^{\frac{n(W)}{\dim W}}$, o\`u $\gg$ signifie ''sup\'erieur \`a, \`a une constante multplicative pr\`es ind\'ependante de $\ell$". En simplifiant par $n(W)$ nous voulons donc montrer que $[L:K]\ll \ell^{\dim W}$. Nous allons en fait montrer mieux : nous allons voir que 
\[\ell^{\dim W}\gg\ll [L:K].\]
\noindent Posons
\[G_H=\left\{ t\in\T(\Z_{\ell})\ | \ \forall i\in I_W\ \ \chi_i(t)=1\mod \ell\right\}.\]
\noindent On a 

\begin{align*}
G_H\cap G_{\ell}	& = \left\{ t\in G_{\ell}\ | \ \forall i\in I_W\ \ \chi_i(t)=1\mod \ell \right\}\\
					& = \left\{ \sigma\in G_{\ell}\ | \ \forall i\in I_W\ \ \sigma(P_i)=P_i\right\}\\
					& = \left\{ \sigma\in G_{\ell}\ | \ \sigma_{|H}=\text{Id} \right\}.
\end{align*}
\noindent Notamment ceci montre que 
\[ [L:K]=\left|\frac{G_{\ell}}{G_{\ell}\cap G_H}\right|.\]
\noindent Or le th\'eor\`eme \ref{trr} dit que l'inclusion $G_{\ell}\hookrightarrow \T(\Z_{\ell})$ est de conoyau de cardinal major\'e par une constante $C_1$ ind\'ependante de $\ell$. On en d\'eduit que l'inclusion de $G_{\ell}/G_{\ell}\cap G_H$ dans $\T(\Z_{\ell})/G_H$ est \'egalement de conoyau major\'e ind\'ependamment de $\ell$. Notamment
\[ [L:K]\gg\ll\left|\frac{\T(\Z_{\ell})}{G_{H}}\right|.\]
\noindent On voudrait maintenant conclure en appliquant le th\'eor\`eme \ref{rib} de Ribet. Pour cela il nous reste \`a interpr\'eter le quotient $\left|\frac{\T(\Z_{\ell})}{G_{H}}\right|$ comme le cardinal de $T_1(\Z/\ell \Z)$ o\`u $T_1$ est un tore alg\'ebrique \`a d\'eterminer. Si $T_1$ est de dimension $\dim W$ on aura alors le r\'esultat attendu~:
\[ [L:K]\gg\ll \ell^{\dim T_1}=\ell^{\dim W}.\]

\noindent Notons $X^*$ le sous-$\Z$-module de $X^*(\T)$ engendr\'e par les $\chi_i$ pour $i\in I_W$ et posons
\[T^{(1)}=X^{*\perp}=\bigcap_{\chi\in X^*}\ker \chi\ \ \text{ et }\ \ T_1=\T/T^{(1)}.\]
\noindent On a une isog\'enie $\T\sim T^{(1)}\times T_1$. Ces groupes alg\'ebriques sont d\'efinis sur $\Q_{\ell}$ par construction, mais \'etant d\'efinis par les caract\`eres $\chi_i$, ils sont \'egalement d\'efinis sur une extension finie de $\Q$. Autrement dit, ils sont d\'efinis sur une extension finie $\kk$ de $\Q$, contenue dans $\Q_{\ell}$. Donc l'isog\'enie est une isog\'enie sur $\kk$. De plus cette extension $\kk/\Q$ est une sous-extension de $\Q(I)$, corps de d\'efinition des $\chi\in I$. Il n'y a donc qu'un nombre fini de telles extensions $\kk$ lorsque $\ell$ varie (parmi les premiers totalement d\'ecompos\'es comme indiqu\'e au d\'ebut du paragraphe).

\medskip

\noindent De plus $X^*(T_1)=X^*\left(\T/X^{*\perp}\right)\simeq X^*$ et $T_1$ est de dimension $\dim W$. Par construction de $T^{(1)}$, on a~:
\[\left|\frac{\T(\Z_{\ell})}{G_H}\right| \gg\ll \left|\frac{T_1(\Z_{\ell})}{\left\{t\in T_1(\Z_{\ell})\ | \ \forall \chi\in X^*\ \chi(t)=1\mod\ell\right\} }\right|.\]
\noindent Suivant les notations de Ribet on a donc
\[\left|\frac{\T(\Z_{\ell})}{G_H}\right| \gg\ll \left|\frac{T_1(\Z_{\ell})}{T_1(1+\ell\Z_{\ell})}\right|\gg\ll \ell^{\dim T_1}=\ell^{\dim W}.\]

\noindent Suivant la m\^eme id\'ee nous allons dans la suite prouver l'in\'egalit\'e $\gamma(A)\leq \alpha(A)$. On commence pour cela par montrer qu'il suffit de savoir montrer l'in\'egalit\'e $\text{Card}\left(H\right)\ll [K(H):K]^{\alpha(A)}$ pour tout sous-groupe fini $H$ de $A[\ell^{\infty}]$, stable par $\Gal(\K/K)$ et pour tout premier $\ell$ : c'est la r\'eduction au cas $\ell$-adique. Ensuite la preuve suit ce qui vient d'\^etre fait pr\'ec\'edemment : on commence par d\'ecouper le groupe $H$ selon les morceaux irr\'eductibles de la repr\'esentation du groupe de Mumford-Tate dans $\GL(\Q_{\ell})$. Il reste ensuite une derni\`ere \'etape consistant \`a recomposer astucieusement ces morceaux de sorte \`a pouvoir estimer le degr\'e de l'extension  des morceaux recompos\'es en fonction de la taille des morceaux. L\`a encore le degr\'e des extensions qui interviendront s'interpr\'etera comme le cardinal des points d'ordre une puissance de $\ell$ d'un tore alg\'ebrique, fabriqu\'e par le m\^eme principe que ci-dessus. Dans la derni\`ere \'etape, il y a deux complications qui se pr\'esentent par rapport au cas trait\'e ci-dessus. D'une part le groupe $H$ n'est plus de type $(\ell,\ldots,\ell)$ : les morceaux intervenant vont donc avoir des poids distincts ce qui complique la partie combinatoire ; d'autre part les caract\`eres n'ont aucune raison d'\^etre d\'efinis dans $\Q_{\ell}$, il faudra donc les grouper par paquets. Afin de rendre plus lisible la preuve nous traiterons pour cette derni\`ere \'etape tout d'abord le cas o\`u les caract\`eres sont d\'efinis sur $\Q_{\ell}$ (\textit{i.e.} $\T_{\ell}$ d\'ecompos\'e sur $\Q_{\ell}$) : il s'agit du paragaphe \ref{casdec}. Le cas g\'en\'eral avec groupement par paquets n'induit essentiellement pas de difficult\'e suppl\'ementaire autre que de notation. La preuve est faite au paragraphe \ref{cg}.

\section{Le th\'eor\`eme \ref{thp}\label{p3} : r\'eduction au cas $\ell$-adique}
\noindent Nous allons maintenant montrer la seconde in\'egalit\'e
\begin{equation}\label{*1}
\gamma(A)\leq\alpha(A).
\end{equation}

\medskip

\noindent Quitte \`a augmenter $K$, on peut supposer (et on le fait) que tous les endomorphismes de $A$ sont d\'efinis sur $K$ et que $A/K$ a bonne r\'eduction en toute place (ceci car d'apr\`es Serre-Tate \cite{serretate} une vari\'et\'e ab\'elienne de type CM a potentiellement bonne r\'eduction en toute place). Ceci va nous permettre de nous ramener comme dans Ribet \cite{ribet} au cas $\ell$-adique. Nous rappelons avant cela un r\'esultat classique concernant les vari\'et\'es ab\'eliennes de type CM~:

\medskip

\begin{lemme}\label{torab} Soit $A/K$ une vari\'et\'e ab\'elienne de type CM sur un corps de nombres $K$. Quitte \`a remplacer $K$ par une extension finie $K'$ ne d\'ependant que de $A$, on a~: pour tout ensemble $S$ de points de torsion de $A(\overline{K})$, l'extension $K(S)/K$ est ab\'elienne.
\end{lemme}
\demo Il suffit de montrer que $A_{\tors}$ engendre une extension ab\'elienne sur $K$. Pour v\'erifier ceci il suffit de traiter le cas d'une vari\'et\'e ab\'elienne simple de type CM et de montrer que pour tout entier $n\geq 1$ et tout premier $\ell\geq 2$, l'extension $K(A[\ell^n])/K$ est ab\'elienne (on prend ensuite le compositum). Ceci revient \`a voir que l'image $G_{\ell}$ de $\rho_{\ell}~:\Gal(\overline{K}/K)\rightarrow \Aut(T_{\ell}(A))\subset \Aut(V_{\ell})$ est commutative. Or ceci est vrai : c'est le corollaire 2 p. 502 de \cite{serretate}.\hfill$\Box$

\medskip

\begin{prop}\label{ladique}Pour d\'emontrer l'in\'egalit\'e (\ref{*1}), il suffit de montrer que : il existe une cons\-tan\-te strictement positive $C(A/K)$ ne d\'ependant que de $A/K$ telle que pour tout nombre premier $\ell$ et tout sous-groupe fini $H$ de $A[\ell^{\infty}]$, stable par $\Gal(\K/K)$, on a
\begin{equation}\label{equa}
\text{Card}\left(H\right)\leq C(A/K)[K(H):K]^{\alpha(A)}.
\end{equation}
\end{prop}
\demo Soit $L/K$ une extension finie. Posons 
\[L'=K\left(A(L)_{\tors}\right),\ H=A(L')_{\tors}\ \text{ et pour tout premier $\ell$, }\ H_{\ell}=A(L')_{\tors}\left[\ell^{\infty}\right].\]
\noindent V\'erifions tout d'abord que $H_{\ell}$ est stable par $\Gal(\K/K)$ : si $\sigma\in\Gal(\K/K)$ et $x\in H_{\ell}$, alors $\sigma(x)\in A[\ell^{\infty}]$. Par ailleurs $A$ \'etant de type CM, les extensions engendr\'ees par des points de torsion sont galoisiennes sur $K$ d'apr\`es le lemme \ref{torab} (quitte \`a avoir au d\'epart remplac\'e, ce que l'on suppose ici, $K$ par une extension finie $K'$ ne d\'ependant que de $A$). Or $K(H_{\ell})\subset L'$, donc $\sigma(x)\in A(L')\cap A[\ell^{\infty}]=H_{\ell}$. On peut donc appliquer l'hypoth\`ese de l'\'enonc\'e aux groupes $H_{\ell}$ pour tout premier $\ell$~:
\[\left| H_{\ell} \right|\ll [K(H_{\ell}):K]^{\alpha(A)}.\]
\noindent En notant $\omega(n)$ le nombre de nombre premiers divisant $n$, il vient 
\begin{equation}\label{in}
\left|A(L)_{\tors}\right|=\left|A(L')_{\tors}\right|=|H|=\prod_{\ell}|H_{\ell}|\leq C(A/K)^{\omega\left(|A(L')_{\tors}|\right)}\prod_{\ell}[K(H_{\ell}):K]^{\alpha(A)}.
\end{equation}
\noindent En effet, sur la d\'efinition de $L'$ on voit imm\'ediatement que $A(L)_{\tors}=A(L')_{\tors}$. Par ailleurs, une estimation classique de $\omega(n)$ est la suivante (\textit{cf.} par exemple \cite{tenenbaum} p. 85 \S\ 5.3)~: $\omega(n)\ll\frac{\log n}{\log\log n}$. L'in\'egalit\'e (\ref{in}) peut donc se r\'e\'ecrire
\[\left|A(L)_{\tors}\right|^{1-\frac{C_1(A/K)}{\log\log \left|A(L)_{\tors}\right|}}\leq \prod_{\ell}[K(H_{\ell}):K]^{\alpha(A)}.\]
\noindent Or $H_{\ell}=A(L')[\ell^{\infty}]$. Si pour tout entier $M\geq 1$, si on note $A(L')[M]$ la partie de $M$-torsion de $A(\overline{K})$ rationnelle sur le corps $L'$, et $d_{L'}(M)=[K(A(L')[M]):K]$, alors par la th\'eorie de Serre-Tate \cite{serretate}, on sait que l'extension $K(A(L')[M])/K$ ne peut \^etre ramifi\'ee qu'en des places au-dessus de premiers divisant $M$. Ainsi si $m$ et $M$ sont premiers entre eux, alors $K(A(L')[m])\cap K(A(L')[M])\subset K'$ o\`u $K'$ est le corps de classes de Hilbert de $K$ (rappelons, cf. par exemple \cite{sil} p.118 Example 3.3, que par d\'efinition le corps de classes de Hilbert de $K$ est la plus grande extension ab\'elienne de $K$ partout non ramifi\'ee). En rempla\c{c}ant $K$ par son corps de classes de Hilbert on obtient ainsi : la fonction
\[M\mapsto d_{L'}(M)\]
\noindent est multiplicative au sens arithm\'etique. On en d\'eduit
\[\left|A(L)_{\tors}\right|^{1-\frac{C_1(A/K)}{\log\log \left|A(L)_{\tors}\right|}}\leq [K(H):K]^{\alpha(A)}\leq [L':K]^{\alpha(A)}\leq [L:K]^{\alpha(A)}.\]
\noindent Ceci conclut.\hfill$\Box$

\medskip

\noindent Dans le paragraphe \ref{dernier} suivant, nous donnons une preuve de l'in\'egalit\'e (\ref{equa}).

\section{Le th\'eor\`eme \ref{thp} : cas de la $\ell^{\infty}$-torsion\label{dernier}}

\subsection{Un r\'esultat pr\'eliminaire}

\noindent Le r\'esultat principal de cette section est la proposition \ref{ld1}.

\subsubsection{Pr\'eliminaires}

\begin{lemme}\label{fl0} Soient $\kk$ un corps, $G/\kk$ un groupe alg\'ebrique, $V$ un $\kk$-espace vectoriel non nul de dimension finie et $\rho : G\rightarrow \GL_V$ une repr\'esentation telle que $\rho : G(\kk)\rightarrow \GL(V)$ est irr\'eductible. Soit $\mathcal{G}$ un sous-groupe abstrait de $G(\kk)$, Zariski dense dans $G$. Alors la repr\'esentation 
\[\rho_{|\mathcal{G}} : \mathcal{G}\rightarrow \GL(V)\]
\noindent est irr\'eductible.
\end{lemme}
\demo Soit $W\not=\{0\}$ un sous-$\kk$-espace vectoriel de $V$ stable par $\rho_{|\mathcal{G}}$. Notons $G_W$ le stabilisateur de $W$ dans $G$ (\textit{i.e.} le foncteur d\'efini par $G_W(R):=\{g\in G(R)\ | \ \rho(g)(W\otimes_{\kk}R)=W\otimes_{\kk}R\}$) : c'est un sous-groupe alg\'ebrique de $G/\kk$. Par hypoth\`ese, $\mathcal{G}\subset G_W(\kk)$, donc par densit\'e de $\mathcal{G}$ dans $G$, on en d\'eduit que $G_W=G$. Ainsi $W$ est stable par la repr\'esentation  $\rho : G(\kk)\rightarrow \GL(V)$, donc par irr\'eductibilit\'e de cette derni\`ere, $W=V$.\hfill $\Box$

\medskip

\noindent Reprenons les notations des paragraphes pr\'ec\'edents. Soit $\ell$ un nombre premier. Le tore $\T/\Q$ op\`ere fid\`ele\-ment sur $V$ par la repr\'esentation $\rho$. On note comme pr\'ec\'edemment $I=\{\chi_1,\ldots,\chi_{2g}\}$ les poids du tore, \textit{i.e.} les caract\`eres diagonalisant l'action de $\T$ sur $V_{\overline{\Q}}$.  

\noindent Si $\chi\in I$ nous noterons $\sigma_1(\chi)=\textnormal{Id},\ldots,\sigma_{t_{\chi}}(\chi)$ les diff\'erents plongements de $\Q_{\ell}(\chi)$, corps de d\'efinition de $\chi$ sur $\Q_{\ell}$, dans $\overline{\Q_{\ell}}$. Sur $\Q_{\ell}$, on peut d\'ecomposer la repr\'esentation 
\[\ \T(\Q_{\ell})\overset{\rho_{\Q_{\ell}}}{\longrightarrow}\GL(V_{\ell}) \]
\noindent en une somme de repr\'esentations irr\'eductibles, chacune se d\'ecomposant \`a nouveau sur $\overline{\Q}_{\ell}$ en $\begin{pmatrix}
\chi	& 		& 				\\
	& \ddots	& 				\\
	& 		& \chi^{\sigma_{t_{\chi}}(\chi)} 	\\
\end{pmatrix}$ o\`u $\chi$ ainsi que ses conjugu\'es sont des \'el\'ements de $I$. Par le lemme \ref{fl0} pr\'ec\'edent (appliqu\'e avec $\kk=\Q_{\ell}$, $G=\T$ et $\mathcal{G}=G_{\ell}$), cette d\'ecomposition en irr\'eductibles correspond \'egalement \`a la d\'ecomposition en irr\'eductibles de la repr\'esentation restreinte
\[\rrho : G_{\ell}\subset \T(\Q_{\ell})\overset{\rho_{\Q_{\ell}}}{\longrightarrow}\GL(V_{\ell}) \]
\noindent correspondant \`a la repr\'esentation $\ell$-adique $\rho_{\ell}$.
 
\medskip

\noindent \'Etant donn\'e un caract\`ere $\chi\in I$, nous noterons $\rho_{\chi}$ (et $V_{\chi}$ le $\Q_{\ell}$-espace vectoriel correspondant) la sous-repr\'esentation irr\'eductible de $\rho_{\Q_{\ell}}$ (ou de mani\`ere \'equivalente $\rrho$) dans laquelle appara\^it $\chi$ quand on \'etend les scalaires \`a $\overline{\Q}_{\ell}$. Le nombre $t_{\chi}$ est la dimension de $V_{\chi}$ sur $\Q_{\ell}$. Les caract\`eres $\chi$ \'etant deux \`a deux distincts, les sous-repr\'esentations $V_{\chi}$ sont uniquement d\'etermin\'ees. On note alors
\[T_{\chi}=V_{\chi}\cap T_{\ell}(A).\]
\noindent Si $x\in V_{\chi}$, on a : $\forall t\in\T(\Q_{\ell})$, $\rho(t)x=\rho_{\chi}(t)x$. De plus la repr\'esentation $\rrho$ (d\'efinie sur $G_{\ell}=\text{Im}(\rho_{\ell})$) est en fait \`a valeur dans $\Aut(T_{\ell}(A))$. En particulier, pour tout $t\in G_{\ell}$, $\rho(t)$ laisse stable $T_{\ell}(A)$ et donc $\rho_{\chi}(t)$ laisse stable $T_{\chi}$.

\medskip

\noindent On a la d\'ecomposition $V_{\ell}=\bigoplus V_{\chi}$, donc en posant $T=\bigoplus T_{\chi}$ on obtient un r\'eseau de $V_{\ell}$ qui est un sous-r\'eseau de $T_{\ell}(A)$ . En prolongeant (pour $t\in \T(\Q_{\ell})$) $\rho_{\chi}(t)$ \`a $V_{\ell}$ par $\rho_{\chi}(t)x=0$ si $x$ est dans $V_{\chi'}$ pour $\chi'$ non conjugu\'e \`a $\chi$, on a 
\[\forall t\in G_{\ell},\ \forall x\in T,\ \ \rho_{\chi}(t)x\in T_{\chi}\subset T.\]

\medskip

\begin{lemme}\label{fl3}L'inclusion $T\subset T_{\ell}(A)$ est de conoyau d'indice born\'e ind\'ependamment de $\ell$. 
\end{lemme}
\demo Rappelons tout d'abord que le module $T_{\ell}(A)$ provient de $\Z$ : en effet on a $T_{\ell}(A)=\H_1(A(\C),\Z)\otimes_{\Z}\Z_{\ell}$ (\textit{cf.} \cite{milne} 15.1 (a)). Les $2g$ caract\`eres $\chi$ sont tous d\'efinis sur une extension finie $\Q(I)$ de $\Q$. Donc la d\'ecomposition de $\rho_{\Q_{\ell}}$ (et donc de $\rrho$) en irr\'eductibles selon les $\rho_{\chi}$ est en fait d\'ej\`a valable sur une sous-extension $L$ de $\Q(I)$, contenue dans $\Q_{\ell}$ ; les $V_{\chi}$ provenant de la d\'ecomposition en irr\'eductibles de $V_L=V\otimes_{\Q} L$ par extension des scalaires \`a $\Q_{\ell}$. Il n'y a qu'un nombre fini de sous-$\Q$-extensions contenues dans $\Q(I)$. Pour $\ell$ variable, les d\'ecompositions de $\rrho$ en irr\'eductibles proviennent donc d'un nombre fini de d\'ecompositions. Ceci montre en particulier que pour $\ell$ assez grand $T=T_{\ell}(A)$. \hfill$\Box$

\medskip

\noindent Nous noterons dans la suite $\ell_0$ un premier tel que si $\ell\geq \ell_0$, alors $T_{\ell}(A)=T$ et nous noterons $n_0$ un entier tel que pour tout premier $\ell$, on a $\ell^{n_0}T_{\ell}(A)\subset T$.

\medskip

\noindent Faisons quelques remarques sur $T$ et les $T_{\chi}$~(les points $1.$ et $2.$ d\'ecoulant de ce que $T=\bigoplus T_{\chi}$)~:
\begin{enumerate}
\item  Soit $m\geq 1$. On a $T_{\chi}/T_{\chi}\cap \ell^mT=T_{\chi}/\ell^mT_{\chi}$.
\item Soit $m\geq 1$. On a $T/\ell^mT=\bigoplus T_{\chi}/\ell^mT_{\chi}.$
\item Soit $n\geq m\geq 1$. On a un isomorphisme
\[\ell^mT_{\chi}/\ell^nT_{\chi}\simeq T_{\chi}/\ell^{n-m}T_{\chi},\]
\noindent donn\'ee par $x\mapsto {\ell}^{-m}x \mod {\ell}^{n-m}$.
\item Enfin on a de m\^eme, pour $n\geq m\geq 1$~:
\[ \left(T_{\chi}/{\ell}^nT_{\chi}\right)[{\ell}^m]={\ell}^{n-m}T_{\chi}/{\ell}^nT_{\chi} \subset  T_{\chi}/{\ell}^{m}T_{\chi}\]
\noindent la derni\`ere inclusion \'etant due au point 3 ; et de m\^eme pour $T$ et pour $T_{\ell}(A)$.
\end{enumerate}

\subsubsection{La proposition}

\noindent Soit $H$ un sous-groupe fini non trivial de $A[\ell^{\infty}]$ stable par Galois (donc par $\rrho$) et d'exposant $\ell^n$. Nous allons d\'ecomposer $H$ selon les sous-repr\'esentations $\rho_{\chi}$. Par d\'efinition $n\geq 1$ est le plus petit entier tel que $H\subset A[\ell^n]$. Si $m\geq 1$, notons 
\[j_m : T_{\ell}(A)/\ell^mT_{\ell}(A)\rightarrow T/\ell^mT, \ \ x\mapsto \begin{cases} \ell^{n_0}x 	& \text{ si }\ \ell<\ell_0\\
																							x			&\text{ si}\ \ell\geq \ell_0.\end{cases}\]
\noindent Notons que si $\ell\geq \ell_0$, l'application $j_m$ n'est autre que l'identit\'e sur $A[\ell^m]$.
\noindent Un premier $\ell$ quelconque \'etant fix\'e, nous allons travailler avec $j_n(H)$ plut\^ot qu'avec $H$. Les remarques pr\'ec\'edentes entra\^inent que $j_n(H)\simeq j_m(H)$ pour $m\geq n$. On renote $H'$ ce groupe. Il est stable par $\rrho$ et $|H'|\gg\ll |H|$, $\gg\ll$ signifiant comparable \`a des constantes ind\'ependantes de $\ell$ et $n$ pr\`es. 

\medskip

\noindent Notons pour tout entier $m\geq 1$, $\pi_m : T\rightarrow T/\ell^mT$ la projection canonique. On d\'efinit
\[\widehat{H'}:=\pi_n^{-1}(H'),\ \ \widehat{H_{\rho_{\chi}}}:=\widehat{H'}\cap V_{\chi},\ \ \text{ et }\ \ H_{\rho_{\chi}}:=\pi_n(\widehat{H_{\rho_{\chi}}}).\]
\noindent Les $H_{\rho_{\chi}}$ sont inclus dans $H'$ et leur somme est directe. De plus, par d\'efinition on a bien le d\'ecoupage voulu : 
\[H_{\rho_{\chi}}=\left\{x\in H'\ |\ \forall t\in G_{\ell},\ \rho(t)x=\rho_{\chi}(t)x\right\}.\]

\medskip

\defi Nous noterons dans la suite $n_{\chi}$ l'entier tel que $\ell^{n_{\chi}}$ est l'exposant de $H_{\rho_{\chi}}$ pour tout $\chi\in I$.

\medskip

\begin{prop}\label{ld1} Le $\Z_{\ell}$-module $\widehat{H'}$ est, \`a un indice fini born\'e ind\'ependamment de $\ell$ pr\`es, la somme directe des $\widehat{H_{\rho_{\chi}}}$ correspondant aux diff\'erentes sous-repr\'esentations irr\'e\-duc\-tibles $\rho_{\chi}$ de $\rrho$. Il en est de m\^eme pour l'inclusion $\bigoplus H_{\rho_{\chi}}\subset H'$. De plus, il existe $\ell_1$ premier et $n_1\in\mathbb{N}^*$ (ind\'ependants de $H$) tels que chaque $H_{\rho_{\chi}}$ est, en tant que groupe, de la forme suivante~:
\begin{equation}\label{e1}
\textit{si }\ell\geq \ell_1\ \textit{ alors }\ H_{\rho_{\chi}}\simeq\left(\Z/\ell^{n_{\chi}}\Z\right)^{t_{\chi}},
\end{equation}
\noindent et,
\begin{equation}\label{e2}
\textit{ si }\ell <\ell_1\ \textit{ et, si }n_{\chi}\geq n_1, \ \textit{ alors }\ H_{\rho_{\chi}}\simeq\prod_{i=1}^{t_{\chi}}\Z/\ell^{n_{\chi}-r_i}\Z
\end{equation}
\noindent o\`u pour tout $i$, $0\leq r_i<n_1$.
\end{prop}

\medskip

\noindent Avant d'aborder la preuve proprement dite, nous commen\c{c}ons par trois lemmes. Les deux premiers seront utilis\'es pour prouver la premi\`ere assertion concernant la d\'ecomposition en somme directe ; le dernier nous servira \`a prouver l'affirmation concernant le cardinal des $H_{\rho_{\chi}}$ pour les petits $\ell$.

\medskip

\noindent On sait par le th\'eor\`eme \ref{trr} qu'il existe une constante $c\in \N^{*}$ telle que 
\[\forall \ell \text{ premier, } \left|\T(\Z_{\ell})/G_{\ell}\right| \text{ divise }c.\]

\medskip

\begin{lemme}\label{fl1} Il existe $t\in \T(\Q)$ tel que $\textnormal{disc}\rho(t^c)\not=0$.
\end{lemme}
\demo Notons $\Delta$ la sous-vari\'et\'e de $\T/\Q$ d\'efinie par $\textnormal{disc} \rho=0$. C'est une sous-vari\'et\'e stricte de  $\T$ (par exemple en passant sur $\overline{\Q}$, on voit que $\Delta_{\overline{\Q}}$ est une r\'eunion finie d'hypersurfaces, car les caract\`eres $\chi$ sont deux \`a deux distincts). Par ailleurs, l'ensemble $\T(\Q)$ est Zariski-dense dans $\T$, donc il en est de m\^eme pour l'ensemble $\{t^c\in\T(\Q)\ |\ t\in \T(\Q)\}$. On peut donc trouver un $t$ comme annonc\'e.\hfill$\Box$

\medskip

\noindent Nous fixons dans la suite $t_0$ un \'el\'ement de $\T(\Q)$ fourni par le lemme \ref{fl1}. 

\medskip

\rem \label{rfl1} Notons que $t_0$ \'etant d\'efini sur $\Q$, il existe un premier $\ell_2$ tel que 
\[\forall \ell\geq \ell_2 \text{ premier }, \ \ \ t_0\in\T(\Z_{\ell}).\]
\noindent En particulier, pour tout premier $\ell\geq \ell_2$ on a  $t_0^c\in G_{\ell}$.

\medskip

\begin{lemme}\label{fl1bis} Le premier $\ell$ \'etant fix\'e, il existe $t_{\ell}\in G_{\ell}$ tel que $\textnormal{disc}\rho(t_{\ell})\not=0$.
\end{lemme}
\demo C'est la m\^eme que pr\'ec\'edemment en travaillant sur $\Q_{\ell}$ plut\^ot que sur $\Q$ et en utilisant la densit\'e de $G_{\ell}$ dans $\T\times_{\Q}\Q_{\ell}$. \hfill$\Box$

\medskip

\begin{lemme}\label{fl2}Soient $\ell$ un nombre premier, $X$ un $\Q_{\ell}$-espace vectoriel de dimension finie $d\geq 1$, $G$ un groupe (abstrait) ab\'elien et $\sigma : G\rightarrow \GL(X)$ une repr\'esentation lin\'eaire irr\'eductible. Soit enfin $R$ un r\'eseau de $X$ (\textit{i.e.} un sous $\Z_{\ell}$-module de $X$ tel que $X=R\otimes_{\Z_{\ell}}\Q_{\ell}$) stable par $\sigma$. Il existe $c_1>0$ (d\'ependant de $\ell$) tel que $\forall n\in\N^*$, $\forall x\in R/\ell^nR$ d'ordre $\ell^n$,  on a~:
\[\left|R/\ell^nR\right|=\ell^{nd}\geq \textnormal{Card}\left(<G\cdot x>\right)\geq c_1\ell^{nd},\]
\noindent o\`u $<G\cdot x>$ est le sous-groupe de $R/\ell^nR$ engendr\'e par les $\sigma(g)(x)$ pour $g\in G$.
\end{lemme}
\demo Pour tout entier $n\geq 1$, notons $\pi_n : R\rightarrow R/\ell^nR$ la r\'eduction modulo $\ell^n$. Supposons par l'absurde le r\'esultat faux : il existe une suite $(x_k)\in R/\ell^{n_k}R$ d'\'el\'ements d'ordre $\ell^{n_k}$ telle que 
\begin{equation}\label{ee1}
\frac{\ell^{n_kd}}{|<G\cdot x_k>|}\underset{k\to +\infty}{\rightarrow}+\infty.
\end{equation}

\noindent Pour tout entier $k$, notons $S_k:=\pi_{n_k}^{-1}(<G\cdot x_{n_k}>)$. La suite $(S_k)$ est une suite de sous-$\Z_{\ell}$-modules de $R$ stables par $\sigma$. \'Ecrivons la d\'ecomposition  de $S_k$ selon les diviseurs \'el\'ementaires~:
\[S_k=\bigoplus_{i=1}^d\ell^{n_i(k)}e_i(k)\Z_{\ell}\]
\noindent avec $n_1(k)\leq\ldots\leq n_d(k)$ et, pour tout $i$, $\ell\nmid e_i(k)$ (si pour un $k_0$ l'un des $e_i(k_0)=0$, alors $S_{k_0}\otimes_{\Z_{\ell}}\Q_{\ell}$ est un sous-espace vectoriel strict non nul de $X$, stable par $\sigma$, ce qui est impossible par irr\'eductibilit\'e : on suppose donc que tous les $e_i(k)$ sont non nuls). S'agissant d'entiers positifs, quitte \`a extraire une sous-suite, on peut supposer (et on le fait) que les suites $(n_i(k))_{k\geq 1}$ sont croissantes pour tout $i$. De plus, comme $x_{n_k}$ est d'ordre $\ell^{n_k}$ dans $R/\ell^{n_k}R$, on voit que $n_1(k)=0$ pour tout $k$. On a
\[|<G\cdot x>|=\textnormal{Card}(\pi_{n_k}(S_k))=\prod_{i=1}^{j(k)}\ell^{n_k-n_i(k)}\]
\noindent o\`u $1\leq j(k)\leq d$ est le plus petit entier tel que $n_{j(k)+1}(k)\geq n_k$ (et o\`u l'on pose $n_{d+1}(k)=+\infty$ pour tout $k$).

\noindent Si la suite $(n_d(k))_{k\geq 1}$ ne tendait pas vers $+\infty$, alors elle serait born\'ee (et donc il en serait de m\^eme des suites $(n_i(k))\leq(n_d(k))$) : 
\[\exists m_0\geq 1\ \forall k\geq 1,\ \ \forall 1\leq i\leq d,\ \ n_i(k)\leq m_0.\]
\noindent Ainsi, si $k$ est assez grand on aurait $n_k > m_0\geq n_d(k)$. Notamment, pour $k$ assez grand on aurait $j(k)=d$ et,
\[|<G\cdot x>|=\ell^{n_kd-\sum_{i=1}^{d}n_i(k)}\geq \ell^{dn_k-dm_0}.\]
\noindent Ceci contredit (\ref{ee1}), donc la suite $(n_d(k))$ tend vers $+\infty$. Par ailleurs, les suite $(n_i(k))$ sont stationnaires ou tendent vers $+\infty$ (et $(n_1)$ est stationnaire et $(n_d)$ tend vers $+\infty$). Notons $a\in\{1,\ldots,d-1\}$ le plus petit entier tel que la suite $(n_{a+1})\to +\infty$. Comme $(n_1)\leq \cdots\leq (n_d)$, on a donc : pour $1\leq i\leq a$, la suite $(n_i)$ est stationnaire de limite not\'ee $n_i$, et, pour $i\geq a+1$, $(n_i)\to+\infty$.

\medskip

\noindent Pour tout entier $k\geq 1$, notons  $\B_k=\{e_1(k),\ldots,e_d(k)\}$. La suite $(\B_k)$ est une suite de bases de $R$. Par compacit\'e on en extrait une sous-suite convergente de limite $\B:=\{e_1,\ldots,e_d\}$. Posons
\[S=\bigoplus_{i=1}^a\ell^{n_i}e_i\Z_{\ell}.\]
\noindent Si $S$ est stable par $\sigma$, alors, comme $S\otimes_{\Z_{\ell}}\Q_{\ell}$ est non nul et de dimension major\'ee par $a$ avec $1\leq a< d$, on pourra conclure par l'absurde gr\^ace \`a l'irr\'eductibilit\'e de $\sigma$. Reste \`a voir que $S$ est stable par $\sigma$. Pour cela, il suffit de voir que pour tout entier $n\geq 1$, la projection $\pi_n(S)$ est stable. Soit donc $n\geq 1$. Nous allons montrer que pour $k$ assez grand,
\[\pi_n(S)=\pi_n(S_k).\]
\noindent Ceci permettra de conclure car $S_k$ est stable par $\sigma$.

\noindent Par stationnarit\'e, il existe $k_0$ tel que $\forall k\geq k_0$, $n_i(k)=n_i$ pour tout $i\in\{1,\ldots,a\}$. De plus il existe $k_0(n)\geq k_0$ tel que $\forall k\geq k_0(n)$, $n_d(k)\geq\ldots\geq n_{a+1}(k)>n$.  Enfin, pour tout $i\in\{1,\ldots,d\}$, la suite $(\pi_n(e_i(k)))_{k\geq 1}$ est une suite convergente d'\'el\'ements de l'ensemble fini $R/\ell^nR$. Elle est donc stationnaire. Donc il existe $k_1(n)\geq k_0(n)$ tel que $\forall k\geq k_1(n)$, $\pi_n(e_i(k))=\pi_n(e_i)$ pour tout $i$. De ceci, on d\'eduit que si $k\geq k_1(n)$, on a $\pi_n(S)=\pi_n(S_k)$. Avec ce qui pr\'ec\`ede, ceci conclut.\hfill$\Box$

\medskip

\noindent \textbf{D\'emonstration de la proposition \ref{ld1}} Commen\c{c}ons par prouver la premi\`ere assertion. La seconde (concernant $\bigoplus H_ {\rho_{\chi}}\subset H'$) s'en d\'eduit imm\'ediatement. Il restera ensuite encore \`a v\'erifier l'in\'egalit\'e num\'erique sur le cardinal des $H_{\rho_{\chi}}$. 

\medskip

\noindent Consid\'erons la repr\'esentation $\rrho : G_{\ell}\subset \T(\Q_{\ell})\overset{\rho_{\Q_{\ell}}}{\rightarrow} \GL(V_{\ell})$. Comme d\'ej\`a expliqu\'e, la d\'ecomposition en irr\'eductibles de $\rrho$ est la m\^eme que celle de $\rho_{\Q_{\ell}}$ (il s'agit du lemme \ref{fl0}).  Notons $V_{\ell}=\bigoplus_{i=1}^r V_{\ell,i}$ cette d\'ecomposition en irr\'eductibles sur $\Q_{\ell}$ et notons $\pr_{\ell,i} : V_{\ell} \rightarrow V_{\ell,i}$ les projecteurs correspondants. Nous noterons \'egalement encore $\pr_{\ell,i}$ les endomorphismes de $V_{\ell}$ obtenus par composition avec l'inclusion $V_{\ell,i}\subset V_{\ell}$. Soit $x\in \widehat{H'}$. On a $x=\sum_{i=1}^r\pr_{\ell,i}(x)$ avec $\pr_{\ell,i}(x)\in V_{\ell,i}$. Par d\'efinition des $\widehat{H_{\rho_{\chi}}}$, il nous suffit donc juste de montrer que $\pr_{\ell,i}(x)$ appartient \`a $\widehat{H'}$ (au moins pourvu que $\ell$ soit assez grand, et, si $\ell$ est petit, que cette d\'ecomposition vaut quitte \`a multiplier au d\'epart $x$ par une constante $M$). Le $\Z_{\ell}$-module $\widehat{H'}$ \'etant stable par $\rrho$, il suffit pour cela de montrer que $\pr_{\ell,i}$ est un polyn\^ome, \`a coefficients dans $\Z_{\ell}$, en $\rho(t)$ pour $t\in G_{\ell}$ convenable. 

\medskip

\noindent On a d\'ej\`a  dit dans la preuve du lemme \ref{fl3} que, $\rho_{\Q_{\ell}}$ provenant de $\rho : \T\rightarrow \GL_V$ par extension des scalaires \`a $\Q_{\ell}$, la d\'ecomposition en irr\'eductible de $\rho_{\Q_{\ell}}$ est d\'ej\`a valable sur une extension finie $L/\Q$ incluse dans $\Q(I)$ ; les $V_{\chi}$ provenant de la d\'ecomposition en irr\'eductibles de $V_L=V\otimes_{\Q} L$ par extension des scalaires \`a $\Q_{\ell}$. Il n'y a qu'un nombre fini de sous-$\Q$-extensions contenues dans $\Q(I)$ (et lorsque $\ell$ varie, les d\'ecompositions de $\rrho$ en irr\'eductibles proviennent donc d'un nombre fini de d\'ecompositions). Notons $V_L=\bigoplus_{i=1}^rV_i$ la d\'ecomposition en irr\'eductibles selon $\rho$, et, pour $i$ compris entre $1$ et $r$, notons $\pr_i : V_L\rightarrow V_i$ les projecteurs correspondants. Nous noterons encore $\pr_i$ les endomorphismes de $V_{L}$ obtenus par composition avec l'inclusion $V_i\subset V_L$. Comme pour les $V_{\ell,i}$, les $\pr_{\ell,i}$ sont d\'efinis \`a partir des $\pr_i$ par extension des scalaires de $L$ \`a $\Q_{\ell}$.

\medskip

\noindent \textbf{Fait : }$\forall t\in \T(L)$,\ \ $\pr_i\circ\rho(t)=\rho(t)\circ\pr_i$.

\noindent En effet, si $y\in V_{L}$, on a
\begin{align*}
\pr_i(\rho(t)y)	& =\pr_i\left(\rho(t)\sum_{j=1}^ry_j\right)= \pr_i\left(\sum_{j=1}^r\rho_j(t)y_j\right)= \rho_i(t)y_i= \rho_i(t)\pr_i(y)=\rho(t)\pr_i(y).
\end{align*}

\medskip

\noindent Rappelons un r\'esultat d'alg\`ebre lin\'eaire : sur un corps $\mathbb{K}$ toute matrice appartenant au commutant d'une matrice donn\'ee $N$ est un polyn\^ome en $N$ \`a coefficients dans $\mathbb{K}$ si et seulement si les polyn\^omes caract\'eristique et minimal de $N$ sont \'egaux. Dans notre situation, pour tout $t$, $\rho(t)$ est diagonalisable sur $\overline{\Q}$. Donc si pour un $t\in \T(L)$ on a $\textnormal{disc}\rho(t)\not=0$ alors les $\pr_i$ sont des polyn\^omes \`a coefficients dans $L$ en $\rho(t)$. Or le lemme \ref{fl1} et la remarque 7.6 nous fournissent un premier $\ell_2$ et un $t_0\in \T(\Q)$ tels que 
\[\text{ si }\ell\geq \ell_2,\ \text{ alors  }\ t_0^c\in G_{\ell},\ \ \text{ et de plus on a } \textnormal{disc}(\rho(t_0^c))\not=0.\]
\noindent Donc $\pr_i$ est un polyn\^ome \`a coefficients dans $L$ en $\rho(t_0^c)$. Notamment il existe un premier $\ell_3\geq\ell_2$ tel que : si $\ell\geq \ell_3$, alors  $\pr_{\ell,i}$ est un polyn\^ome \`a coefficients $\Z_{\ell}$ en $\rho(t_0^c)$ (le premier $\ell_3$ d\'epend a priori de $\pr_i$ donc du $\ell$ fix\'e au d\'ebut de la section, mais comme il n'y a qu'un nombre fini de d\'ecompositions possibles lorsque $\ell$ varie, on peut en fait prendre pour $\ell_3$ une valeur ind\'ependante de $\ell$). Comme $t_0^c\in G_{\ell}$ ceci conclut : si $\ell\geq \ell_3$, alors $\widehat{H'}$ est la somme directe des $\widehat{H_{\rho_{\chi}}}$. Si $\ell<\ell_3$, alors on applique le lemme \ref{fl1bis} en lieu et place du lemme \ref{fl1} (en travaillant directement sur $\Q_{\ell}$ plut\^ot que sur $L$) et on obtient que $\pr_{\ell,i}$ est un polyn\^ome \`a coefficients dans  $\Q_{\ell}$ en $\rho(t_{\ell})$. Il existe donc une constante $M\in\Z$ telle que $M\pr_{\ell,i}\in \Z_{\ell}[\rho(t_{\ell})]$ (l\`a encore $M$ est a priori d\'ependante de $\ell$ mais comme $\ell<\ell_3$ on peut \'evidemment choisir $M$ ind\'ependante de $\ell$). Donc $M\pr_{\ell,i}(x)\in\widehat{H}$ et donc $M\widehat{H}\subset \bigoplus \widehat{H_{\rho_{\chi}}}$. Ceci prouve donc la premi\`ere assertion.

\medskip

\noindent D\'eterminons maintenant la forme de $H_{\rho_{\chi}}$. En tant que groupe, $H_{\rho_{\chi}}$ est un produit de $\Z/\ell^{m}\Z$ comportant au plus $t_{\chi}$ facteurs et avec des $m\leq n_{\chi}$. L'assertion est donc une question de cardinalit\'e.

\medskip

\noindent Rappelons que si $\ell\geq \ell_0$, alors $T=T_{\ell}(A)$. Il existe $\ell_1\geq \ell_0$ premier tel que pour tout premier $\ell\geq \ell_1$, la d\'ecomposition $\rrho\simeq\bigoplus_{i=1}^m\rho_i$ reste irr\'eductible sur $\F_{\ell}$. En effet, soit $\rho_{\chi}$ l'une des sous-repr\'esentations irr\'eductibles sur $\Q_{\ell}$ consid\'er\'ees. Elle se diagonalise sur $\overline{\Q}_{\ell}$ selon l'ensemble des caract\`eres $E=\{\chi^{\sigma}\ |\ \sigma\in\Gal(\overline{\Q}_{\ell}/\Q_{\ell})\}$. La repr\'esentation $\rho_{\chi}\mod\ell$ se diagonalise sur $\overline{\F}_{\ell}$ selon $\overline{E}=\{\chi^{\sigma}\ |\ \sigma\in\Gal(\overline{\F}_{\ell}/\F_{\ell})\}$. Dire que $\rho_{\chi}\mod\ell$ reste irr\'eductible sur $\F_{\ell}$ correspond donc \`a dire que $\overline{E}$ ne se d\'ecompose pas en plusieurs orbites sous l'action de $\Gal(\overline{\F}_{\ell}/\F_{\ell})$. Comme $E$ forme une unique orbite on en d\'eduit qu'il en est de m\^eme pour $\overline{E}$ pour tout $\ell$ suffisamment grand. Ceci va nous permettre de d\'eterminer la forme de $H_{\rho_{\chi}}$. 

\medskip

\noindent Commen\c{c}ons par le cas des $\ell<\ell_1$ et fixons $\chi$. Le lemme \ref{fl2} appliqu\'e avec $G=G_{\ell}$, $X=V_{\chi}$ de dimension $t_{\chi}$, $R=T_{\chi}$ et $\sigma=\rho_{\chi} : G_{\ell}\rightarrow GL(V_{\chi})$ donne le r\'esultat annonc\'e (en prenant un point $x\in H_{\rho_{\chi}}$ d'ordre maximal $\ell^{n_{\chi}}$).

\medskip

\noindent Passons au cas $\ell\geq \ell_1$. En r\'eduisant modulo $\ell H'$ on voit, par irr\'eductibilit\'e et cardinalit\'e, que si $H_{\rho_{\chi}}$ n'est pas le groupe trivial $\{0\}$, alors sa r\'eduction modulo $\ell$, que nous noterons ici plus simplement $H_{1}$, est isomorphe \`a $\prod_{i=1}^{t_{\chi}}\Z/\ell\Z$. Il reste \`a voir ce qu'il advient quand on passe de $H_{1}$ \`a $H_{\rho_{\chi}}$. Le groupe $H_{\rho_{\chi}}$ contient un \'el\'ement d'ordre $\ell^{n_{\chi}}$. De plus $[\ell^{n_{\chi}-1}] H_{\rho_{\chi}}$ est contenu dans $H_{1}$ et est stable par $\rrho$, donc par irr\'eductibilit\'e on voit que $[\ell^{n_{\chi}-1}]~: H_{\rho_{\chi}}\rightarrow H_{1}$ est surjective. Ceci implique que $H_{\rho_{\chi}}$ est de la forme attendue $\prod_{i=1}^{t_{\chi}}\Z/\ell^{n_{\chi}}\Z.$\hfill$\Box$ 
 
\medskip

\noindent \textbf{Notations } Rappelons que dans la suite on note $n_{\chi}$ l'entier tel que $\ell^{n_{\chi}}$ est l'exposant de $H_{\rho_{\chi}}$. Par ailleurs, si $\ell<\ell_1$ et $n_{\chi}\geq n_1$, on notera $n_{\chi}^{\min}=\min n_{\chi}-r_i$ et pour unifier les notations on notera \'eventuellement $n_{\chi}^{\min}:=n_{\chi}$ lorsque $\ell\geq \ell_1$.

\medskip

\noindent Nous voulons dans la suite montrer une in\'egalit\'e du type 
\[|H|\ll[K(H):K]^{\alpha(A)}.\]
\noindent Comme $|H|\gg\ll |H'|$ et comme $[K(H'):K]\leq [K(H):K]$, on voit qu'il suffit de prouver la m\^eme assertion avec $H'$ plut\^ot que $H$. Pr\'ecis\'ement nous allons obtenir une majoration du cardinal de  $\bigoplus H_{\rho_{\chi}}$ en fonction de $[K(\bigoplus H_{\rho_{\chi}}):K]\leq [K(H'):K]\leq [K(H):K]$. Ceci se traduira donc en la m\^eme borne, \`a constante multiplicative pr\`es, pour le cardinal de $H$.

\subsection{Interm\`ede : cas d\'ecompos\'e\label{casdec}}

\noindent Nous donnons ici la preuve pour les premiers $\ell$ tels que les caract\`eres $\chi\in I$ sont d\'efinis sur $\Q_{\ell}$, \textit{i.e.} tels que le tore $\T$ est scind\'e sur $\Q_{\ell}$. Nous supposerons \'egalement $\ell\geq \ell_1$ o\`u $\ell_1$ est la valeur introduite dans la proposition \ref{ld1} pr\'ec\'edent. Bien qu'inutile du strict point de vue logique la preuve dans cette situation particuli\`ere permet d'illustrer le raisonnement dans un cadre plus simple, notamment du point de vue des notations. Consid\'erons donc un tel $\ell$. La repr\'esentation $\rho$ se d\'ecompose sur $\Q_{\ell}$ en $2g$ sous-repr\'esentations irr\'eductibles de dimension $1$ donn\'ees par les caract\`eres $\chi\in I$. Autrement dit on a $\rho_{\chi}=\chi$. Notamment \textit{nous noterons (dans cette section uniquement) $H_{\chi}$ le groupe not\'e $H_{\rho_{\chi}}$ pr\'ec\'edemment}. On a l'isomorphisme de groupes donn\'e par la proposition \ref{ld1}
\[H_{\chi}\simeq \Z/\ell^{n_{\chi}}\Z.\]

\defi Nous dirons que le caract\`ere $\chi\in I$ \textit{intervient dans $H$} si le sous-groupe $H_{\chi}$ de $H$ correspondant \`a $\chi$ est non trivial, \textit{i.e.} si $n_{\chi}\not= 0$.

\medskip

\noindent On pose $W_H=\text{Vect}_{\Q}(\chi_1,\ldots,\chi_s)$ l'espace vectoriel engendr\'e par tous les caract\`eres intervenant dans $H$. \`A chaque tel caract\`ere $\chi$ est associ\'e l'entier $n_{\chi}\geq 1$. On note 
\[ n^{(1)}=\sup_{\chi\in I}n_{\chi},\ \ W_1=\text{Vect}_{\Q}\left(\chi^{(1)}\right), \]
\noindent o\`u $\chi^{(1)}$ est un caract\`ere intervenant dans $H$ associ\'e \`a $n^{(1)}$. De plus on note
\[ I_1=I\cap W_1,\ \  b_1=\text{Card}\left(I_1\right),\ \ H_1=\bigoplus_{\chi\in I_1}H_{\chi}.\]

\medskip

\noindent On d\'efinit alors par r\'ecurrence (jusqu'\`a un rang $r$ tel que $W_r=W_H$) pour tout entier $i\geq 2$
\[ n^{(i)}=\!\!\!\!\sup_{\chi\in I\backslash I_{i-1}}\!\!\!\!n_{\chi},\ \ W_{i}=\text{Vect}_{\Q}\left(W_{i-1},\chi^{(i)}\right),\ \ I_{i}=I\cap W_i,\]
\noindent o\`u $\chi^{(i)}$ est un caract\`ere appartenant \`a $I\backslash I_{i-1}$, associ\'e \`a l'entier $n^{(i)}$. Par ailleurs, on note
\[\forall i\geq 2,\ \ b_{i}=\text{Card}I_{i}-\text{Card}I_{i-1},\ \ H_{i}=\bigoplus_{\chi\in I_{i}}H_{\chi}.\]
\noindent On note $r$ le plus petit entier $i$ tel que $W_i=W_H$.

\medskip

\noindent Par construction, les $W_i$ sont de dimension $i$ et la famille $\{\chi^{(1)},\ldots,\chi^{(r)}\}$ est une base de $W_H$. De plus, pour tout $i$ compris entre $1$ et $r$, le cardinal de $H_i$ est major\'e par $\ell^{\sum_{k=1}^ib_kn^{(k)}}$. Notamment, par la proposition \ref{ld1}, on a 
\begin{equation}\label{pg1dec}
|H|\ll  \ell^{\sum_{k=1}^rb_kn^{(k)}}.
\end{equation}

\noindent Enfin pour tout $1\leq i\leq r$, $\sum_{k=1}^ib_k$ est le nombre de caract\`eres de $I$ appartenant \`a $W_i$. Il nous reste \`a estimer le degr\'e de l'extension $K(H)/K$. Pr\'ecis\'ement nous voulons montrer que 
\begin{equation}\label{pg2dec}
\ell^{\sum_{i=1}^r n^{(i)}}\ll [K(H):K].
\end{equation}

\medskip

\noindent Admettons pour l'instant cette in\'egalit\'e (\ref{pg2dec}) et expliquons comment conclure dans ce cas. Les in\'egalit\'es (\ref{pg1dec}) et (\ref{pg2dec}) impliquent~:
\begin{equation}\label{concluredec}
|H|\ll [K(H):K]^{\frac{\sum_{i=1}^rn^{(i)}b_i}{\sum_{i=1}^rn^{(i)}}}.
\end{equation}

\noindent On utilise alors le lemme suivant :

\medskip

\begin{lemme}\label{lg4}Soient $n_1\geq \ldots\geq n_r$, $b_1,\ldots,b_r$  et $w_1,\ldots,w_r$ des entiers strictement positifs. On a 
\[\frac{\sum_{i=1}^rn_ib_i}{\sum_{i=1}^rn_iw_i}\leq \sup_{1\leq k\leq r}\frac{\sum_{i=1}^kb_i}{\sum_{i=1}^kw_i}.\]
\end{lemme}
\demo On pose $n_{r+1}=0$ et on applique une transformation d'Abel \`a la somme $\sum n_ib_i$ : 
\begin{align*}
\sum_{i=1}^{r+1}n_ib_i	& =\left(\sum_{i=1}^{r+1} b_i\right)n_{r+1}+\sum_{i=1}^r\left(\sum_{k=1}^ib_k\right)(n_i-n_{i+1})= \sum_{i=1}^r\left(\sum_{k=1}^ib_k\right)(n_i-n_{i+1})\\
			& \leq \left(\sup_{1\leq i\leq r}\frac{\sum_{k=1}^ib_k}{\sum_{k=1}^iw_k}\right)\sum_{i=1}^r\left(\sum_{k=1}^iw_k\right)(n_i-n_{i+1}) \leq \left(\sup_{1\leq i\leq r}\frac{\sum_{k=1}^ib_k}{\sum_{k=1}^iw_k}\right)\sum_{i=1}^rw_in_i.\\
\end{align*}
\noindent On conclut en divisant par le nombre strictement positif $\sum_{i=1}^rw_in_i$.\hfill$\Box$

\medskip

\noindent On constate donc finalement que l'exposant maximal obtenu dans l'in\'egalit\'e (\ref{concluredec}) est atteint quand les $n^{(i)}$ valent tous $0$ ou $1$, \textit{i.e.} quand $H$ est en fait un sous-groupe de $A[\ell]$. De plus dans ce cas, un majorant de l'exposant maximal est
\[\sup_{1\leq i\leq r} \frac{\text{nombre de caract\`eres dans }W_i}{\text{dim }W_i}.\]
\noindent Ainsi, cet exposant est inf\'erieur au nombre $\alpha(A)$ pr\'ec\'edemment d\'efini, ce qui conclut.

\medskip

\noindent Il nous reste maintenant \`a prouver l'in\'egalit\'e (\ref{pg2dec}). On reprend pour cela ce qui avait \'et\'e fait dans la section \ref{derp}. Consid\'erons pour cela la tour d'extensions
\[
\xymatrix{
K \ar[r]& \cdots		\ar[r]	& K(H_{1}) 	\ar[r]	& \cdots		\ar[r]	 & K(H_r)		\ar[r]	& K(A[\ell^{\infty}]) }
\]

\noindent Pour tout entier $i$ entre $1$ et $r$, notons $\mathcal{G}_i:=\left\{\sigma\in G_{\ell}\ | \ \sigma_{|H_i}=\textnormal{Id}\right\}$ le groupe de Galois correspondant \`a l'extension $K(H_i)$. On a 
\begin{align*}
\mathcal{G}_i	& =\left\{\sigma\in G_{\ell}\ | \ \forall \chi\in I_i\ \ \forall P_{\chi}\in H_{\chi}\ \ \ \sigma P_{\chi}=P_{\chi}\right\}	\\
				& =\left\{\sigma\in G_{\ell}\ | \ \forall \chi\in I_i\ \ \ \chi(\sigma)=1 \mod \ell^{n_{\chi}}\right\}							\\
				& \subset \left\{\sigma\in G_{\ell}\ | \ \forall k\leq i \ \ \ \chi^{(k)}(\sigma)=1\mod\ell^{n^{(k)}}\right\}						\\
				& =  \left\{t\in \T(\Z_{\ell})\ | \ \forall k\leq i\ \ \  \chi^{(k)}(t)=1\mod\ell^{n^{(k)}}\right\}\cap G_{\ell}					
\end{align*}

\noindent On pose $G_i:=\left\{t\in \T(\Z_{\ell})\ | \ \forall k\leq i\ \ \  \chi^{(k)}(t)=1\mod\ell^{n^{(k)}}\right\}$ de sorte que $\mathcal{G}_i\subset G_i\cap G_{\ell}$. Ainsi
\[[K(H_i):K]=\left|\frac{G_{\ell}}{\mathcal{G}_i}\right|\geq\left|\frac{G_{\ell}}{G_i\cap G_{\ell}}\right|\gg\ll\left|\frac{\T(\Z_{\ell})}{G_i}\right|.\]
\noindent Comme pr\'ec\'edemment nous allons maintenant interpr\'eter ce dernier terme comme le cardinal des points d'un certain tore. On va montrer par r\'ecurrence sur $i\leq r$ que 
\begin{equation}\label{findec}
\left|\frac{\T(\Z_{\ell})}{G_i}\right|\gg \ell^{\sum_{k=1}^in^{(k)}},
\end{equation}
\noindent le cas $i=r$ \'etant la conclusion attendue. Commen\c{c}ons par le cas $i=1$~: pour cela on pose $T^{(1)}=\ker\chi^{(1)}$ et $T_1=\T/T^{(1)}$ . Le groupe des caract\`eres de $T_1$ est engendr\'e par $\chi_1$ et on a
\begin{align*}\left|\frac{\T(\Z_{\ell})}{G_1}\right| 	& \gg\ll\left|\frac{T_1(\Z_{\ell})}{\left\{t\in T_1(\Z_{\ell})\ |\ \chi^{(1)}(t)=1\mod \ell^{n^{(1)}}\right\}}\right|\\
												& \gg\ll \ell^{n^{(1)}}\hspace{1cm} \text{ par le th\'eor\`eme \ref{rib} de Ribet.}
\end{align*}
\noindent Ceci prouve (\ref{findec}) quand $i=1$.

\medskip

\noindent Supposons maintenant l'in\'egalit\'e vraie au rang $i$ et montrons-la au rang $i+1\leq r$. On a
\[\left|\frac{\T(\Z_{\ell})}{G_{i+1}}\right|=\left|\frac{\T(\Z_{\ell})}{G_{i}}\right|\times \left|\frac{G_i}{G_{i+1}}\right|.\]
\noindent En appliquant l'hypoth\`ese de r\'ecurrence il suffit donc de montrer que $\left|\frac{G_i}{G_{i+1}}\right|\gg \ell^{n^{(i+1)}}$. On introduit pour cela
\[\forall i\geq 1\ \ T^{(i+1)}=\ker{\chi^{(i+1)}}_{|T^{(i)}}\ \ \text{ et }\ \ T^{(i)}\sim T^{(i+1)}\times T_{i+1}.\]
\noindent o\`u l\`a encore le groupe des caract\`eres de $T_{i+1}$ est engendr\'e par $\chi^{(i+1)}$. Comme on l'avait not\'e pr\'ec\'edemment ces groupes alg\'ebriques sont d\'efinis sur $\Q_{\ell}$ par construction, mais \'etant d\'efinis par des caract\`eres $\chi\in I$, ils sont \'egalement d\'efinis sur une extension finie de $\Q$. Autrement dit, ils sont d\'efinis sur une extension finie $\kk$ de $\Q$, contenue dans $\Q_{\ell}$. Donc les diff\'erentes isog\'enies sont des isog\'enies sur $\kk$. De plus cette extension $\kk/\Q$ est une sous-extension de $\Q(I)$, corps de d\'efinition des $\chi\in I$. Il n'y a donc qu'un nombre fini de telles extensions $\kk$ lorsque $\ell$ varie (parmi les premiers totalement d\'ecompos\'es).

\medskip

\noindent Par construction de $T^{(i)}$, on a
\[\frac{T^{(i)}(\Z_{\ell})}{T^{(i)}(\Z_{\ell})\cap G_{i+1}}\hookrightarrow \frac{G_i}{G_{i+1}}.\]
\noindent On en d\'eduit
\[\left|\frac{G_i}{G_{i+1}}\right|\geq \left|\frac{T^{(i)}(\Z_{\ell})}{T^{(i)}(\Z_{\ell})\cap G_{i+1}}\right|\gg \left|\frac{T_{i+1}(\Z_{\ell})}{\left\{ t\in T_{i+1}(\Z_{\ell})\ | \ \chi^{(i+1)}(t)=1\mod l^{n^{(i+1)}}\right\}}\right|\gg\ll \ell^{n^{(i+1)}}.\]
\noindent Ceci conclut la preuve dans le cas d\'ecompos\'e. Il nous reste maintenant \`a traiter le cas g\'en\'eral : la preuve est exactement la m\^eme \`a ceci pr\`es que nous devons regrouper les caract\`eres $\chi$ par paquets avec les $\chi^{\sigma}$, $\sigma$ d\'ecrivant les plongements de $\Q_{\ell}(\chi)$ dans $\overline{\Q_{\ell}}$. De m\^eme les tores \`a consid\'erer seront cette fois d\'efinis par $\bigcap_{\sigma}\ker \chi^{\sigma}$ au lieu de $\ker\chi$. Nous aurons notamment besoin d'utiliser le lemme \ref{lg4} avec des $w_i$ non n\'ecessairement tous \'egaux \`a $1$.

\subsection{Cas g\'en\'eral\label{cg}}

\subsubsection{Les objets combinatoires\label{combi}}

\noindent Les caract\`eres $\chi_1,\ldots,\chi_{2g}$ sont d\'efinis sur $\overline{\Q}$ donc sur une extension $L/\Q$ finie.

\medskip

\defi On dit que le caract\`ere $\chi\in I$ \textit{intervient dans $H$} si~ : 
\begin{enumerate}
\item soit $\ell\geq \ell_1$ et le sous-groupe $H_{\rho_{\chi}}$ de $H$ correspondant \`a la repr\'esentation $\rho_{\chi}$ est non trivial.
\item soit $\ell<\ell_1$  et le sous-groupe $H_{\rho_{\chi}}$ de $H$ correspondant \`a la repr\'esentation $\rho_{\chi}$ est d'exposant $\ell^{n_{\chi}}\geq \ell^{n_1}$.
\end{enumerate}

\medskip

\noindent On pose $W_H=\text{Vect}_{\Q}(\chi_1,\ldots,\chi_s)$ l'espace vectoriel engendr\'e par tous les caract\`eres intervenant dans $H$. Ici et dans la suite on suppose $H$ de cardinal suffisamment grand de sorte que $W_H\not=\{0\}$ (si $H$ est petit il n'y a rien \`a montrer). \`A chaque caract\`ere $\chi$ intervenant dans $H$ est associ\'e l'entier $n_{\chi}$. Comme au paragraphe \ref{casdec} pr\'ec\'edent, mais cette fois en regroupant les caract\`eres conjugu\'es, on note 
\[n^{(1)}=\sup_{\chi\in I}n_{\chi},\ \ W_1=\text{Vect}_{\Q}\left((\chi^{(1)})^{\ss_1(1)},\ldots,(\chi^{(1)})^{\ss_{t_1}(1)}\right), \]
\noindent o\`u $\chi^{(1)}$ est un caract\`ere appartenant \`a $I$ associ\'e \`a $n^{(1)}$ et o\`u $\ss_1(1),\ldots,\ss_{t_1}(1)$ sont les diff\'erents plongements de $\Q_{\ell}(\chi^{(1)})$ dans $\overline{\Q}_{\ell}$. De plus on note
\[ I_1=I\cap W_1,\ \ w_1=\dim W_1,\ \  b_1=\text{Card}\left(I_1\right),\ \ H_1=\bigoplus_{\chi\in I_1}H_{\rho_\chi},\]
\noindent o\`u la somme porte sur les $\chi\in I_1$ modulo l'action de Galois. On extrait de l'ensemble $I_1$ une base $\left\{(\chi^{(1)})^{\ss_1(1)},\ldots,(\chi^{(1)})^{\ss_{w_1}(1)}\right\}$ de $W_1$.

\medskip

\rem Notons que si $\chi\in I_1$, alors pour tout plongement $\sigma$, on a $\chi^{\ss}\in I_1$, l'espace $W_1$ \'etant par construction stable par l'action de Galois. 

\medskip

\noindent On d\'efinit alors par r\'ecurrence (jusqu'au rang $r$ tel que $W_r=W_H$) pour tout entier $i\geq 2$
\[ n^{(i)}=\!\!\!\!\sup_{\chi\in I\backslash I_{i-1}}\!\!\!\!n_{\chi},\ \ W_{i}=\text{Vect}_{\Q}\left(W_{i-1},(\chi^{(i)})^{\ss_1(i)},\ldots,(\chi^{(i)})^{\ss_{t_i}(i)}\right),\ \ I_{i}=I\cap W_i,\]
\noindent o\`u $\chi^{(i)}$ est un caract\`ere appartenant \`a $I\backslash I_{i-1}$, associ\'e \`a l'entier $n^{(i)}$ et o\`u $\ss_1(i),\ldots,\ss_{t_i}(i)$ sont les diff\'erents plongements de $\Q_{\ell}(\chi^{(i)})$ dans $\overline{\Q}_{\ell}$. Par ailleurs, on note
\[\forall i\geq 2,\ \  w_i=\dim W_i-\dim W_{i-1},\ \ b_{i}=\text{Card}I_{i}-\text{Card}I_{i-1},\ \ H_{i}=\bigoplus_{\chi\in I_{i}}H_{\rho_\chi},\]
\noindent o\`u  la somme porte sur les $\chi\in I_i$ modulo l'action de Galois. On compl\`ete, avec des \'el\'ements de $I_i$, la base de $W_{i-1}$ (construite par la r\'ecurrence) en une base de $W_i$. On note $r$ le plus petit entier $i$ tel que $W_i=W_H$.

\medskip

\noindent Finalement quitte \`a r\'eordonner les termes, et pour soulager les notations, on suppose que 
\[\forall r\geq i\geq 1,\ \ \sigma_1(i)=\text{Id}, \ \ \chi_i=\chi^{(i)},\  \text{ et donc }\ n_i:=n_{\chi_{i}}=n^{(i)}.\]

\begin{lemme}\label{lg2} La famille $\{\chi_1,\ldots,(\chi_1)^{\ss_{w_1}(1)},\ldots,\chi_r,\ldots,(\chi_r)^{\ss_{w_r}(r)}\}$ est une base de $W_H$. De plus, pour tout $i$ compris entre $1$ et $r$, le cardinal de $H_i$ est major\'e par $\ell^{\sum_{k=1}^ib_kn_k}$ et $\sum_{k=1}^ib_k$ est le nombre de caract\`eres (\'el\'ements de $I$) appartenant \`a $W_i$.
\end{lemme}
\demo Ceci d\'ecoule de la construction des diff\'erents objets et de la proposition \ref{ld1}.\hfill$\Box$

\subsubsection{Degr\'e du corps de rationalit\'e}

\noindent Nous voulons maintenant minorer le degr\'e de l'extension $K(H)/K$. Pr\'ecis\'ement nous voulons prouver l'in\'egalit\'e suivante~:
\begin{equation}\label{pg2} 
\ell^{\sum_{i=1}^rw_in_i}\ll [K(H):K].
\end{equation}

\noindent La conclusion suit alors exactement comme au paragraphe \ref{casdec} de cette in\'egalit\'e et du lemme \ref{lg2} pr\'ec\'edent. En effet, admettons pour l'instant l'in\'egalit\'e (\ref{pg2}). On a donc
\begin{equation}\label{pg1}
\text{Card}H\ll \ell^{\sum_{i=1}^rn_ib_i}\ll [K(H):K]^{\frac{\sum_{i=1}^rn_ib_i}{\sum_{i=1}^rn_iw_i}}.
\end{equation}

\noindent On utilise alors le lemme \ref{lg4} comme on l'avait fait dans le cas d\'ecompos\'e, la seule diff\'erence \'etant que les $w_i$ ne sont plus n\'ecessairement tous \'egaux \`a $1$. On constate finalement que l'exposant maximal obtenu dans l'in\'egalit\'e (\ref{pg1}) est atteint quand les $n_i$ valent tous $0$ ou $1$, \textit{i.e.} quand $H$ est en fait un sous-groupe de $A[\ell]$. De plus dans ce cas, un majorant de l'exposant maximal est
\[\sup_{1\leq i\leq r} \frac{\text{nombre de caract\`eres dans }W_i}{\text{dim }W_i}.\]
\noindent Ainsi, cet exposant est inf\'erieur au nombre $\alpha(A)$ pr\'ec\'edemment d\'efini, ce qui nous permet de conclure la partie $\ell$-adique et donc la preuve de l'in\'egalit\'e $\gamma(A)\leq \alpha(A)$ du th\'eor\`eme \ref{thp}.\hfill$\Box$

\medskip

\noindent Il nous reste finalement \`a prouver l'in\'egalit\'e (\ref{pg2}). Consid\'erons pour cela la tour d'extensions
\[
\xymatrix{
K \ar[r]& \cdots		\ar[r]	& K(H_{1}) 	\ar[r]	& \cdots		\ar[r]	 & K(H_r)		\ar[r]	& K(A[\ell^{\infty}]) }
\]

\noindent Pour tout entier $i$ entre $1$ et $r$, notons $\mathcal{G}_i:=\left\{\sigma\in G_{\ell}\ | \ \sigma_{|H_i}=\textnormal{Id}\right\}$ le groupe de Galois correspondant \`a l'extension $K(H_i)$. 

\medskip

\begin{lemme}\label{galois} On a l'inclusion 
\[\mathcal{G}_i\subset G_{\ell}\cap \bigcap_{k=1}^i\left\{t\in \T(\Z_{\ell})\, | \, \forall j\in [\![1,t_k]\!] \ \chi^{\ss_j(k)}_k(t)=1\mod \ell^{n^{\min}_k}\right\}=:G_{\ell}\cap G_i.\]
\end{lemme}
\demo En utilisant les notations du paragraphe pr\'ec\'edent, on a
\begin{align*}
\mathcal{G}_i & =\left\{t\in G_{\ell}\, | \, \forall \chi\in I_i,\ \ \forall P\in H_{\rho_{\chi}},\ \  \rho_{\chi}(t)P=P\right\}	\\
    & \subset \left\{t\in G_{\ell}\, | \, \forall 1\leq k\leq i,\ \ \  \forall P_k\in H_{\rho_{\chi_k}},\ \  \rho_{\chi_k}(t)P_k=P_k\right\}		\\
    & \subset \left\{t\in G_{\ell}\, | \, \forall 1\leq k\leq i,\ \  \rho_{\chi_k}(t)=1\mod \ell^{n^{\min}_k}\right\}\\
    & =G_{\ell}\cap\bigcap_{k=1}^i\left\{t\in \T(\Z_{\ell})\, | \, \forall j\in [\![1,t_k]\!],\ \ \chi^{\ss_j(k)}_k(t)=1\mod \ell^{n^{\min}_k}\right\}.	\\
\end{align*}
\noindent La derni\`ere \'egalit\'e d\'ecoule du fait que par construction, la repr\'esentation $\rho_{\chi_k}$ est \'equivalente \`a $\begin{pmatrix}
\chi_k	& 		& 				\\
	& \ddots	& 				\\
	& 		& \chi_k^{\ss_{t_k}(k)} 	\\
\end{pmatrix}$ sur $\overline{\Q}_{\ell}$. Ceci conclut.\hfill$\Box$ 

\medskip

\noindent Ainsi
\[[K(H_i):K]=\left|\frac{G_{\ell}}{\mathcal{G}_i}\right|\geq\left|\frac{G_{\ell}}{G_i\cap G_{\ell}}\right|\gg\ll\left|\frac{\T(\Z_{\ell})}{G_i}\right|.\]
\noindent Comme pr\'ec\'edemment nous allons maintenant interpr\'eter ce dernier terme comme le cardinal des points d'un certain tore. Pr\'ecis\'ement nous allons montrer par r\'ecurrence que 
\begin{equation}\label{fin}
\forall 1\leq i\leq r,\ \ \left|\frac{\T(\Z_{\ell})}{G_i}\right|\gg \ell^{\sum_{k=1}^in^{\min}_kw_k}.
\end{equation}
\noindent L'in\'egalit\'e (\ref{pg2}) en d\'ecoule alors avec $i=r$ (en fait il reste encore \`a voir que l'on peut remplacer les entier $n^{\min}_i$ par les $n_i$ ; mais ceci est imm\'ediat, car on sait qu'il existe deux constantes absolues $\ell_1$, $n_1$ telles que pour tout $\ell\geq \ell_1$, $n^{\min}_i=n_i$, et que pour les $\ell < \ell_1$ en nombre fini restants, on a $n^{\min}_i\geq n_i-n_1$) . 

\medskip

\noindent Pour montrer l'in\'egalit\'e (\ref{fin}) consid\'erons pour tout entier $i$ compris entre $1$ et $r$,
\[T^{(0)}:=\T, \ \ \text{ et }\ \ \forall i\geq 0\ \ T^{(i+1)}=\bigcap_{k=1}^{t_{i+1}}\ker \left(\chi^{\ss_k(i+1)}_{i+1}\right)_{|T^{(i)}}\text{ et }\ \ T^{(i)}\sim T^{(i+1)}\times T_{i+1}.\]
\noindent o\`u pour $i\geq 1$ le groupe des caract\`eres de $T_{i}$ est $X_i=(T^{(i)})^{\perp}=\left\{\chi\in X^*(T^{(i-1)})\, | \, T^{(i)}\subset \ker\chi\right\}.$ 

\medskip

\noindent Ces tores sont d\'efinis sur $\Q_{\ell}$ par construction, mais \'etant d\'efinis par les caract\`eres $\chi_i$, ils sont \'egalement d\'efinis sur une extension finie de $\Q$. Autrement dit, ces tores sont d\'efinis sur une extension finie $\kk$ de $\Q$, contenue dans $\Q_{\ell}$. De plus cette extension est une sous-extension de $L$, corps de d\'efinition des $\chi\in I$. Il n'y a donc qu'un nombre fini de telles extensions $\kk$ lorsque $\ell$ varie.

\medskip

\begin{lemme}\label{lc1}Pour tout $i$ compris entre $1$ et $r$ le tore $T_i$ est de dimension $w_i$.
\end{lemme}
\demo Le groupe des caract\`eres de $T_i$ est $X_i$ et par construction, $X_i\otimes \Q$ est engendr\'e par les caract\`eres $\chi_{i}^{\ss_k(i)}$ pour $k$ variant entre $1$ et $t_i$. La dimension de $T_i$ est donc $\dim W_i-\dim W_{i-1}=w_i$.\hfill$\Box$

\medskip

\noindent Nous pouvons maintenant prouver l'in\'egalit\'e (\ref{fin}) par r\'ecurrence sur $i\leq r$. V\'erifions tout d'abord la propri\'et\'e au rang $i=1$~: par construction de $T_1$ et comme dans le cas d\'ecompos\'e, on a 
\[\left|\frac{\T(\Z_{\ell})}{G_1}\right|\gg\ll \left|\frac{T_1(\Z_{\ell})}{T_1(1+\ell^{n^{\min}_1}\Z_{\ell})}\right|.\]
\noindent Le th\'eor\`eme \ref{rib} permet de conclure :
\[\ell^{w_1n^{\min}_1}\ll\left|T_1\left(\Z/\ell^{n^{\min}_1}\Z\right)\right|.\]
\noindent Le passage du rang $i$ au rang $i+1$ est exactement le m\^eme que dans le cas d\'ecompos\'e trait\'e au paragraphe \ref{casdec}. Notons que les diff\'erentes constantes intervenant dans les in\'egalit\'es $\gg$ d\'ependent du corps $\kk$, mais on a vu que lorsque $\ell$ varie, seul un nombre fini de tels $\kk$ interviennent, donc on peut bien choisir dans les $\gg$ des constantes multiplicatives ind\'ependantes de $\ell$.

\section{Les cas particuliers}

\noindent Nous allons traiter deux types de cas particuliers : certains cas particuliers de vari\'et\'es ab\'e\-lien\-nes de type CM d'une part et le cas des courbes elliptiques sans multiplication complexe d'autre part.

\subsection{Cas particuliers pour des vari\'et\'es ab\'eliennes de type CM}

\noindent On rappelle en suivant Kubota \cite{kubota}, qu'une vari\'et\'e ab\'elienne de type CM est de type non d\'eg\'en\'er\'e si $d=g+1$, avec $g=\text{dim }A$ et $d=\text{dim }\T$ o\`u $\T$ est le groupe de Mumford-Tate de $A$.

\medskip

\defi Avec les notations pr\'ec\'edentes en suivant Kubota \cite{kubota}, on dit qu'une vari\'et\'e ab\'elienne de type CM est de \textit{d\'efaut} $\delta$ si 
\[\delta=g+1-d.\]

\noindent Comme pr\'ec\'edemment, on note $\chi_1,\ldots,\chi_{2g}$ les caract\`eres diagonalisant l'action de $\T$ sur $V_{\overline{\Q}}$. Quitte \`a les renum\'eroter, on peut regrouper ces caract\`eres par deux de sorte que
\[\chi_1+\chi_{g+1}=\cdots=\chi_{g}+\chi_{2g}=:\chi_0.\]
\noindent Ceci d\'ecoule par exemple de \cite{del} I Example 3.7 point (d) p. 47. De fait Deligne donne l'\'enonc\'e dual pour le groupe des cocaract\`eres, mais la traduction est imm\'ediate. Il faut \'egalement faire attention au fait que sa d\'efinition du groupe de Mumford-Tate est tr\`es l\'eg\`erement diff\'erente de la n\^otre, ce qui fait appara\^itre un cocaract\`ere $e_0$ qui n'a pas lieu d'\^etre avec nos conventions.

\medskip

\noindent Dans le cas non d\'eg\'en\'er\'e, ces relations sont les seules relations de d\'ependance lin\'eaire. Par ailleurs, \'etant donn\'e un sous-$\Q$-espace vectoriel $W$ non nul de $X^*(\T)\otimes\Q$, on introduit deux notations
\[I_W=\left\{i\in\{1,\ldots,g\}\, | \, \chi_i\in W \text{ ou }\chi_{g+i}\in W\right\},\ \text{ et }\ m(W)=\left|I_W\right|.\]

\begin{lemme}\label{d0}Si $A/K$ est de type non d\'eg\'en\'er\'e, alors
\[\alpha(A)=\frac{2g}{d}.\]
\end{lemme}
\demo On distingue deux cas :
\begin{enumerate}
\item Si $\chi_0\notin W$ alors $n(W)=m(W)$ o\`u l'on rappelle que $n(W)$ est le nombre de caract\`eres, parmi les $\chi_1,\ldots,\chi_{2g}$, appartenant \`a $W$. De plus $\dim W\geq m(W)$, donc dans ce cas on a m\^eme $\frac{n(W)}{\dim W}\leq 1\leq \frac{2g}{d}$.
\item Sinon $\chi_0\in W$. Dans ce cas, on voit que $n(W)=2m(W)$ et $\dim W\geq m(W)+1$. En utilisant la croissance de la fonction $x\mapsto \frac{2x}{x+1}$, on a donc 
\[\frac{n(W)}{\dim W}\leq\frac{2m(W)}{m(W)+1}\leq \frac{2g}{g+1}=\frac{2g}{d}.\]
\end{enumerate}
\noindent Ceci conclut.\hfill$\Box$

\medskip

\noindent De m\^eme qu'il est facile de traiter le cas des vari\'et\'es ab\'eliennes de type non d\'eg\'en\'er\'e, il est facile de traiter le cas des vari\'et\'es ab\'eliennes de d\'efaut $\delta=1$.

\medskip

\begin{lemme}\label{d1}Si $A/K$ est de d\'efaut $1$, alors
\[\alpha(A)=2=\frac{2g}{d}.\]
\end{lemme}
\demo Notons que si $W=X^*(\T)\otimes \Q$, alors $\frac{n(W)}{\dim W}=\frac{2g}{g}=2$. Soit maintenant $W$ un sous-espace vectoriel non nul strict de $X^*(\T)\otimes\Q$. Comme dans le lemme pr\'ec\'edent, on distingue deux cas.
\begin{enumerate}
\item Si $\chi_0\notin W$ alors $n(W)=m(W)$ et $\dim W\geq m(W)-1$. Le quotient est donc inf\'erieur \`a $2$.
\item Sinon $\chi_0\in W$. Dans ce cas $n(W)=2m(W)$ et $\dim W\geq m(W)$. L\`a encore, le quotient est inf\'erieur \`a $2$, ce qui conclut.\hfill$\Box$
\end{enumerate}

\medskip

\noindent Nous pouvons maintenant passer \`a la preuve de la proposition \ref{parti1} : notons tout d'abord qu'il suffit en fait de montrer la majoration de $\alpha(A)$. L'\'egalit\'e d\'ecoule alors de notre th\'eor\`eme \ref{thp} et de l'in\'egalit\'e (\ref{tribcor}). Ceci \'etant, si la dimension $g$ est un nombre premier ou $1$, alors le type est non d\'eg\'en\'er\'e et donc le lemme \ref{d0} permet de conclure. Il reste \`a traiter le cas des vari\'et\'es ab\'eliennes simples de dimension $4$ ou $6$.

\medskip

\begin{enumerate}
\item Si $g=4$, on sait par la remarque \ref{cds} que la dimension $d$ du groupe de Mumford-Tate est comprise entre $2+\log_2(4)=4$ et $4+1=5$. Si $d=5$ le type est non d\'eg\'en\'er\'e et si $d=4$ la vari\'et\'e est de d\'efaut $\delta=1$. L\`a encore les deux lemmes pr\'ec\'edents permettent de conclure.
\item Si $g=6$, l'encadrement pr\'ec\'edent vaut toujours : $d$ est compris entre $4<2+\log_2(6)$ et $6+1=7$. Comme pr\'ec\'edemment toujours, si $g=7$ ou si $g=6$, les lemmes \ref{d0} et \ref{d1} pr\'ec\'edents permettent de conclure. Il reste le cas o\`u $d=5$. Dans ce cas, par d\'efinition $\alpha(A)=\sup \frac{n(W)}{\dim W}$ o\`u le sup porte sur les sous-espaces vectoriels non nuls de $X^*(\T)\otimes \Q$. En prenant pour $W$ l'espace tout entier on constate que $\alpha(A)\geq \frac{2g}{d}\geq 2$. Par ailleurs, si $W$ est un sous-espace strict de $X^*(\T)\otimes\Q$, sa dimension est inf\'erieure \`a $4$. Le corollaire \ref{c1new} du paragraphe \ref{p2} nous indique que $n(W)\leq 2^{4-1}=8$. Ainsi on a $\frac{n(W)}{\dim W}\leq \frac{8}{4}=2\leq \frac{2g}{d}$. Ceci prouve bien que $\alpha(A)=\frac{2g}{d}$ et ach\`eve la preuve de la proposition \ref{parti1}.\hfill$\Box$
\end{enumerate}

\subsection{Cas d'une courbe elliptique sans multiplication complexe}

\noindent Soit $E/K$ une courbe elliptique sans multiplication complexe. Comme pr\'ec\'edemment, on peut, en utilisant l'ind\'ependance alg\'ebrique des repr\'esentations $\ell$-adiques $\rho_{\ell}$ (\textit{cf.} \cite{serrenoncm} th\'eor\`eme 3), se ramener au cas $\ell$-adique : on se donne un groupe fini $H$ inclus dans $E[\ell^{\infty}]$ pour un certain premier $\ell$. La courbe $E/K$ \'etant sans multiplication complexe, on sait par un r\'esultat de Serre (th\'eor\`eme 3 et son corollaire 1 de \cite{serrenoncm} p. 299--300) que pour presque tout premier $\ell$, on a 
\[G_{\ell}=\GL_2(\Z_{\ell}).\]
\noindent Ainsi au vu de ce que l'on veut montrer, on peut supposer que $G_{\ell}=\GL_2(\Z_{\ell})$. Le groupe $H$ est de la forme
\[H=\langle P_1\rangle\oplus\langle P_2\rangle\simeq \Z/\ell^{n_1}\Z\times\Z/\ell^{n_2}\Z.\]
\noindent Par ailleurs, le groupe de Galois associ\'e, $G_H$, tel que $[G_{\ell}:G_H]=[K(H):K]$ est d\'efini par 
\[G_H=\left\{\sigma\in G_{\ell}\ / \ \sigma_H=\text{Id}_H\right\}.\]
\noindent On se donne une base de $T_{\ell}(E)$, $\{\hat{P_1},\hat{P_2}\}$ telle que $P_i=\hat{P_i}\mod \ell^{n_i}$ pour $i\in\{1,2\}$. On a ainsi l'identification
\[G_H=\left\{\sigma\in \GL_2(\Z_{\ell})\, | \, \forall i\in\{1,2\}\ \ \sigma \hat{P_i}=\hat{P_i}\mod \ell^{n_i}\right\}.\]
\noindent On peut encore r\'e\'ecrire ceci sous la forme
\[G_H=\left\{ \begin{pmatrix}
a	& b	\\
c	& d	\\
\end{pmatrix}\in \GL_2(\Z_{\ell})\ / \ a-1=c=0\mod \ell^{n_1}\ \text{ et }\ b=d-1=0\mod \ell^{n_2}\right\}.\]
\noindent Sur cette derni\`ere \'ecriture, on voit qu'il existe deux constantes absolues $C_1$ et $C_2$ strictement positives telles que 
\[C_1\leq\frac{[G_{\ell}:G_H]}{\ell^{2(n_1+n_2)}}\leq C_2.\]
\noindent Le groupe $H$ \'etant pr\'ecis\'ement de cardinal $\ell^{n_1+n_2}$, ceci permet de conclure.\hfill$\Box$

\vspace{1cm}

\noindent \textbf{Adresse :} Ratazzi Nicolas, \\ Universit\'e Paris Sud\\ D\'epartement de math\'ematiques, B\^atiment 425\\ 91405 Orsay Cedex \\ FRANCE


\begin{thebibliography}{10}

\bibitem{BMZ}
E.~Bombieri, D.~Masser et U.~Zannier.
\newblock Intersecting a curve with algebraic
  subgroups of multiplicative groups.
\newblock {\em Internat. Math. Res. Notices}, volume~20, pages 1119--1140,
  1999.

\bibitem{bmz2}
E.~Bombieri, D.~Masser et U.~Zannier.
\newblock Intersecting curves and algebraic subgroups : conjectures and more results.
\newblock {\em Trans. Amer. Math. Soc}, volume~358, pages 2247--2257,
2006.

\bibitem{boro}
M.~V. Borovo{\u\i}.
\newblock The action of the {G}alois group on the rational cohomology classes
  of type {$(p,\,p)$} of abelian varieties.
\newblock {\em Mat. Sb. (N.S.)}, 94(136):649--652, 656, 1974.

\bibitem{dodson}
B.~Dodson.
\newblock On the {M}umford-{T}ate group of an abelian variety with complex
  multiplication.
\newblock {\em J. Algebra}, 111(1):49--73, 1987.

\bibitem{del}
P. Deligne, J.~S. Milne, A. Ogus et K.-Y. Shih.
\newblock {\em Hodge cycles, motives, and {S}himura varieties}, volume 900 of
  {\em Lecture Notes in Mathematics}.
\newblock Springer-Verlag, Berlin, 1982.

\bibitem{kubota}
T.~Kubota.
\newblock On the field extension by complex multiplication.
\newblock {\em Trans. Amer. Math. Soc.}, 118:113--122, 1965.

\bibitem{lettre}
D.~Masser.
\newblock Lettre \`a \textnormal{Daniel Bertrand} du 10 novembre 1986.

\bibitem{mas}
D.~Masser.
\newblock Small values of the quadratic part of the
  \textnormal{N}\'eron-\textnormal{T}ate height.
\newblock In {\em Progr. Math.}, volume~12, pages 213--222.
  Birkh\text{\"a}user, 1981.

\bibitem{mazur}
B.~Mazur.
\newblock Rational isogenies of prime degree (with an appendix by {D}.
  {G}oldfeld).
\newblock {\em Invent. Math.}, 44(2):129--162, 1978.

\bibitem{merel}
L.~Merel.
\newblock Bornes pour la torsion des courbes elliptiques sur les corps de
  nombres.
\newblock {\em Invent. Math.}, 124(1-3):437--449, 1996.

\bibitem{milne}
J.~S. Milne.
\newblock Abelian varieties.
\newblock In {\em Arithmetic geometry (Storrs, Conn., 1984)}, pages 103--150.
  Springer, New York, 1986.

\bibitem{moonen}
B.~Moonen et Y.~Zarhin.
\newblock Hodges classes on abelian varieties of low dimension.
\newblock {\em Math. Ann.}, 315, n. 4:711--733, 1999.

\bibitem{mumf}
D.~Mumford.
\newblock A note of {S}himura's paper ``{D}iscontinuous groups and abelian
  varieties''.
\newblock {\em Math. Ann.}, 181:345--351, 1969.

\bibitem{murty}
V.K. Murty.
\newblock Hodge and {W}eil classes on abelian varieties.
\newblock In {\em The arithmetic and geometry of algebraic cycles (Banff, AB,
  1998)}, volume 548 of {\em NATO Sci. Ser. C Math. Phys. Sci.}, pages 83--115.
  Kluwer Acad. Publ., Dordrecht, 2000.

\bibitem{ono}
T.~Ono.
\newblock Arithmetic of algebraic tori.
\newblock {\em Ann. Math.}, 74, n. 1:101--139, 1961.

\bibitem{parent}
P.~Parent.
\newblock Bornes effectives pour la torsion des courbes elliptiques sur les
  corps de nombres.
\newblock {\em J. Reine Angew. Math.}, 506:85--116, 1999.

\bibitem{pink1}
R.~Pink.
\newblock $\ell$-adic algebraic monodromy groups cocharacters, and the {M}umford-{T}ate conjecture.
\newblock {\em J. Reine Angew. Math.}, 495:187--237, 1998. 

\bibitem{pink}
R.~Pink.
\newblock A common generaliza\-tion of the conjectures of {A}ndr\'e-{O}ort,
  {M}anin-{M}umford, and {M}ordell-{L}ang. 
\newblock Pr\'epublication de 2005 dis\-po\-nible \`a l'adres\-se
  http://www.math.ethz.ch/~pink/ftp/AOMMML.pdf.

\bibitem{piat}
I.~I. Pjatecki{\u\i}-{\v{S}}apiro.
\newblock Interrelations between the {T}ate and {H}odge hypotheses for abelian
  varieties.
\newblock {\em Mat. Sb. (N.S.)}, 85(127):610--620, 1971.

\bibitem{pohlmann}
H.~Pohlmann.
\newblock Algebraic cycles on abelian varieties of complex multiplication type.
\newblock {\em Ann. Math.}, 88, n. 2:161--180, 1968.

\bibitem{mray}
M. Raynaud.
\newblock Courbes sur une vari\'et\'e ab\'elienne et points de torsion.
\newblock {\em Inv. Math}, 71:207--233, 1983.

\bibitem{semilehmer}
N.~Ratazzi.
\newblock Minoration de la hauteur sur les vari\'et\'es ab\'eliennes de type
  \textnormal{CM} et applications.
\newblock Pr\'e\-pub\-li\-ca\-tion de l'{I}nstitut de math\'ematiques de {J}ussieu,
  f\'evrier 2005.

\bibitem{remond}
G.~R\'emond.
\newblock Intersection de sous-groupes et de sous-vari\'et\'es \textnormal{I}.
\newblock {\em Math. Ann.}, 333, n. 3:525--548, 2005.

\bibitem{rv}
G.~R\'emond et E.~Viada.
\newblock Probl\`eme de {M}ordell-{L}ang modulo certaines sous-vari\'et\'es ab\'eliennes,
\newblock {\em Int. Math. Res. Not.}, 35:1915--1931, 2003.

\bibitem{ribet}
K.~A. Ribet.
\newblock Division fields of abelian varieties with complex multiplication.
\newblock {\em M\'emoires de la S.M.F}, 2:75--94, 1980.

\bibitem{ribet2}
K.~A. Ribet.
\newblock Hodge classes on certain types of abelian varieties.
\newblock {\em Amer. J. Math.}, 105:523--538, 1983.

\bibitem{serrenoncm}
J.-P. Serre.
\newblock Propri\'et\'es galoisiennes des points d'ordre fini des courbes
  elliptiques.
\newblock {\em Inv. Math}, 15:259--331, 1972.


\bibitem{serladique}
J.-P. Serre.
\newblock Repr\'esentations $\ell$-adiques.
\newblock In {\em Kyoto Int. Symposium on Algebraic Number Theory}, pages
  177--193. Japan Soc. for the promotion of Science, 1977.

\bibitem{serken}
J.-P. Serre.
\newblock {\em Lettre \`a Ken Ribet, {\OE}uvres. {C}ollected papers. {IV}}.
\newblock Springer-Verlag, Berlin, 2000.
\newblock 1985--1998.

\bibitem{college}
J.-P. Serre.
\newblock {\em R\'esum\'e des cours au coll\`ege de France de 1984-1985, {\OE}uvres. {C}ollected papers. {IV}}.
\newblock Springer-Verlag, Berlin, 2000.
\newblock 1985--1998.

\bibitem{serabl}
J.-P. Serre.
\newblock {\em Abelian {$\ell$}-adic representations and elliptic curves},
volume~7 of {\em Research Notes in Mathematics},
\newblock With the collaboration of Willem Kuyk and John Labute, Revised reprint of the 1968 original,
\newblock A K Peters Ltd., Wellesley, MA, 1998.


\bibitem{serretate}
J.-P. Serre et J.~Tate.
\newblock Good reduction of abelian varieties.
\newblock {\em Ann. Math.}, 88:492--517, 1968.

\bibitem{sga4}
{\em Th\'eorie des topos et cohomologie \'etale des sch\'emas. {T}ome 3}.
\newblock Springer-Verlag, Berlin, 1973.
\newblock S\'eminaire de G\'eom\'etrie Alg\'ebrique du Bois-Marie 1963--1964
  (SGA 4), Dirig\'e par M. Artin, A. Grothendieck et J. L. Verdier. Avec la
  collaboration de P. Deligne et B. Saint-Donat, Lecture Notes in Mathematics,
  Vol. 305.

\bibitem{shita}
G.~Shimura et Y.~Taniyama.
\newblock {\em Complex multiplication of abelian varieties and its applications
  to number theory}, volume~6 of {\em Publications of the Mathematical Society
  of Japan}.
\newblock The Mathematical Society of Japan, Tokyo, 1961.

\bibitem{silv}
A.~Silverberg.
\newblock Torsion points on abelian varieties of {CM}-type.
\newblock {\em Compositio Math.}, 68(3):241--249, 1988.

\bibitem{sil}
J.~H.~Silverman.
\newblock {\em Advanced Topics in the Arithmetic of Elliptic Curves}, 
\newblock Graduate Texts in Mathematics, Vol.~156, 1999.

\bibitem{tenenbaum}
G.~Tenenbaum.
\newblock {\em Introduction \`a la th\'eorie analytique et probabiliste des
  nombres}, volume~1 de {\em Cours Sp\'ecialis\'es}.
\newblock Soci\'et\'e Math\'ematique de France, Paris, seconde \'edition, 1995.

\bibitem{viada}
E.~Viada.
\newblock The intersection of a curve with algebraic subgroups in a product of
  elliptic curves.
\newblock In {\em Ann. Scuola Norm. Pisa Cl. Sci. \textnormal{S\'erie} (V)},
  volume~2, pages 47--75, 2003.

\bibitem{japo}
H.~Yanai.
\newblock On the rank of {CM}-type.
\newblock {\em Nagoya Math. J.}, 97:169--172, 1985.

\bibitem{zilber}
B.~Zilber.
\newblock Exponential sums equations and the {S}chanuel conjecture.
\newblock {\em J. London Math. Soc. (2)}, 65(1):27--44, 2002.


\end{thebibliography}
\end{document}